\documentclass{amsproc}
\usepackage{amsfonts,amsmath,amstext,mathrsfs,epsfig,color}
      
\newtheorem{theorem}{Theorem}[section]
\newtheorem{lemma}[theorem]{Lemma}
\newtheorem{proposition}[theorem]{Proposition}
\newtheorem{corollary}[theorem]{Corollary}
\theoremstyle{definition}
\newtheorem{definition}[theorem]{Definition}
\theoremstyle{remark}
\newtheorem{remark}[theorem]{Remark}
\newtheorem{example}[theorem]{Example}

\renewcommand{\enspace}{\;}
\newcommand{\sP}{\mathcal{P}}
\newcommand{\NEW}[1]{{\em #1}\index{#1}}
\newcommand{\mrm}[1]{\text{#1}}
\newcommand{\argmax}{\operatornamewithlimits{arg\, max}}       
\def\<#1,#2>{\langle #1, #2\rangle}

\def\intt#1{\mathop{\mathrm{int}} (#1)}
\def\domain{\mathop{\mathrm{dom}}}
\def\lsc{\mathrm{lsc}}
\def\usc{\mathrm{usc}}
\def\supp{\mathop{\mathrm{idom}}}

\def\impl{\mathop{\;\Rightarrow\;}}

\def\dom#1{\domain (#1)}
\def\idom#1{\supp (#1)}
\def\ridom#1{\mathop{\mathrm{ridom}} (#1)}
\def\ldom#1{\mathop{\mathrm{ldom}}(#1)}
\def\udom#1{\mathop{\mathrm{udom}}(#1)}

\def\sF{{\mathscr F}}
\def\sFC{{\mathscr F}_{\!c}}

\def\sG{{\mathscr G}}
\def\AC{{\mathscr A}}
\def\cP{{\mathscr P}}
\def\sS{{\mathscr S}}
\def\sB{{\mathscr B}}
\def\pX{\cP(X)}
\def\pZ{\cP(Z)}
\def\pW{\cP(W)}
\def\pY{\cP(Y)}
\def\R{\mathbb{R}}
\def\N{\mathbb{N}}

\newcommand{\leqf}{\leq_{\sF}}
\newcommand{\leqg}{\leq_{\sG}}
\newcommand{\gal}{^{\circ}}
\newcommand{\bgal}{b\gal}
\newcommand{\set}[2]{\{#1\mid\,#2\}}

\newcommand{\minf}{\text{\rm min}_{\sF}}
\newcommand{\ming}{\text{\rm min}_{\sG}}
\newcommand{\supf}{\text{\rm sup}_{\sF}}
\newcommand{\infg}{\text{\rm inf}_{\sG}}
\newcommand{\inff}{\text{\rm inf}_{\sF}}
\newcommand{\supg}{\text{\rm sup}_{\sG}}
\newcommand{\rmaxb}{\overline{\R}_{\max}}

\newcommand{\rbar}{\overline{\R}}
\newcommand{\bydef}{\stackrel{\text{\rm def}}{=}}
\newcommand{\RBY}{\rbar^{_{\scriptstyle Y}}}
\newcommand{\RBZ}{\rbar^{_{\scriptstyle Z}}}
\newcommand{\RBBZ}{\rmaxb^{_{\scriptstyle Z}}}
\newcommand{\RBX}{\rbar^{_{\scriptstyle X}}}
\newcommand{\RBYP}{\rbar^{_{\scriptstyle \{y_1,y_2\}}}}
\newcommand{\RBYPj}{\rbar^{_{\scriptstyle \{y_j\}}}}
\newcommand{\condc}{$({\mathcal{C}})$}
\newcommand{\condcmath}{({\mathcal{C}})}

\newcommand{\lipo}[1]{\mathrm{Lip}_{\omega}(#1)}
\newcommand{\sgn}{\operatorname{sgn}}
\newcommand{\myb}{\sgn(\lambda)\lambda^2}
\newcommand{\mybinv}{\sgn(\lambda)\sqrt{|\lambda|}}
\newcommand{\topsf}{\top_{\!\sF}}
\newcommand{\topS}{\top_{\!S}}

\makeindex
\begin{document}
\title[Invertibility of Functional Galois Connections]{Set coverings and invertibility of Functional Galois Connections}
\author{Marianne Akian}
\address{Marianne Akian,
INRIA, Domaine de Voluceau,
B.P.~105, 78153 Le Chesnay Cedex, France.}
\email{Marianne.Akian@inria.fr}
\author{St\'ephane Gaubert}
\address{St\'ephane Gaubert,
INRIA, Domaine de Voluceau,
B.P.~105, 78153 Le Chesnay Cedex, France.}
\email{Stephane.Gaubert@inria.fr}
\author{Vassili Kolokoltsov}
\address{Vassili Kolokoltsov,
Dep.\ of Computing and Mathematics,
Nottingham Trent University, Burton Street,
Nottingham, NG1 4BU, UK, 
and Institute for Information Transmission Problems of Russian
Academy of Science, Moscow, Russia.}
\email{vk@maths.ntu.ac.uk}
\thanks{{\em Date}: July 26, 2003. {\em ESI--preprint}:
January 28, 2004. {\em Revised}: March 25, 2004}
\thanks{This work has been partially supported by the
Erwin Schr\"odinger International Institute of Mathematical Physics (ESI)} 
\begin{abstract}
We consider equations of the form $Bf=g$,
where $B$ is a Galois connection between
lattices of functions. This includes the case
where $B$ is the Legendre-Fenchel transform,
or more generally a Moreau conjugacy.
We characterise the existence and uniqueness
of a solution $f$ in terms of 
generalised subdifferentials.
This extends a theorem of Vorobyev and Zimmermann,
relating solutions of max-plus linear equations
and set coverings.
We give various illustrations.
\end{abstract}
\maketitle

\section{Introduction}
We call (dual) {\em functional Galois connection}
a (dual) Galois connection between a sublattice
$\sF$ of $\RBY$ and a sublattice $\sG$ of $\RBX$,
where $X,Y$ are two sets and $\rbar=\R\cup\{\pm\infty\}$,
see Section~\ref{sec-rep} for definitions.
An important example of functional Galois
connection is the Legendre-Fenchel transform, and
more generally, the Moreau
conjugacy~\cite{moreau70}
associated to a kernel $b: X\times Y \to \rbar$,
\begin{align}
B: \RBY\to\RBX,
\qquad Bf(x)= \sup\set{b(x,y)-f(y)}{y\in Y}\enspace,
\label{moreau}
\end{align}
where $b(x,y)-f(y)$ is an abbreviation of $b(x,y)+(-f(y))$,
with the convention that $-\infty$ is absorbing for addition,
i.e., $-\infty+\lambda=\lambda+(-\infty)=-\infty$,
for all $\lambda\in\rbar$.
Moreau conjugacies are instrumental
in nonconvex duality, see~\cite[Chapter~11, Section~E]{rockwets},\cite{singer}.
Max-plus linear operators with kernel,
which are of the form $f\mapsto B(-f)$,
arise in deterministic optimal control and asymptotics,
and have been widely studied, see in 
particular~\cite{cuning,maslov92,BCOQ,maslovkolokoltsov95,DENSITE,gondran}.  
Other examples of functional Galois connections include 
dualities for quasi-convex functions (see for 
instance~\cite{singer,volle,singer2}).

We consider here general (dual) Galois connections between the 
set $\sF$ of lower semicontinuous functions from a Hausdorff topological
space $Y$ to $\rbar$ and 
$\sG=\RBX$. Any functional Galois
connection $B:\sF\to\sG$ has the form:
\[ Bf (x)=\sup \set{b(x,y,f(y))}{y\in Y}\enspace,\]
where $b: X\times Y\times \rbar \to \rbar$ is such that
$b(x,\cdot,\alpha)\in \sF$, for all $x\in X$ and $\alpha\in\rbar$,
and $b(x,y,\cdot)$ is nonincreasing, right continuous and sends
$+\infty$ to $-\infty$, for all $x\in X$ and $y\in Y$,
see Theorem~\ref{theo-rep} and Proposition~\ref{prop-blsc} below.
Representation theorems of this type have been
obtained previously by Maslov and 
Kolokoltsov~\cite{KOLO88,KOLO90,kolokoltsov92} and 
Mart\'\i nez-Legaz and Singer~\cite{martinez90,singer},
see Section~\ref{sec-rep}.

Given a map $g\in \sG$ and 
a functional Galois connection $B:\sF\to\sG$,
we consider the problem:
\[
(\sP):\;\;\; 
\mrm{Find }f\in \sF\mrm{ such that } Bf =g\enspace.
\]
In particular, we look for effective conditions on $g$ for 
the solution $f$ to exist and be unique. 
When $X,Y$ are finite sets, $\sF=\R^Y$, $\sG=\R^X$, 
and $B$ is as in~\eqref{moreau} with a real valued kernel $b$,
the solutions of $(\sP)$ were characterised by 
Vorobyev~\cite[Theorem  2.6]{vorobyev67} in terms of ``minimal 
resolvent coverings'' of $X$ (in fact, Vorobyev 
considered equations of the type $\min_{y\in Y} a(x,y) f(y)=g(x)$ 
where $a,f$ and $g$ take (finite) nonnegative values, which correspond
in~\eqref{moreau}, to kernels $b$ with values in $\R\cup\{+\infty\}$). 
This approach was systematically 
developed by Zimmermann~\cite[Chapter 3]{Zimmermann.K}, who 
considered several algebraic structures and 
allowed in particular the kernel $b$ to take the value  $-\infty$.
The method of Vorobyev and Zimmermann
is one of the basic tools in max-plus linear algebra,
and it has been instrumental in the understanding of the 
geometry of images of max-plus linear operators.
It has been the source of several important 
developments, including the characterisation 
by Butkovi\v{c}~\cite{butkovip94,butkovic} of locally
injective (``strongly regular'') max-plus linear maps in terms
of optimal assignment problems.

We extend here Vorobyev's theorem
to the case of functional Galois connections.
We use an adapted notion of subdifferentials, which is similar 
to that introduced by Mart\'\i nez-Legaz and Singer~\cite{martinez95}.
In the special case of Moreau conjugacies,
subdifferentials were introduced by Balder~\cite{balder},
Dolecki and Kurcyusz~\cite{dolecki}, and Lindberg~\cite{lindberg}
(see also~\cite{martinez88}).
We show that the existence (Section~\ref{exis-sec})
and uniqueness  (Section~\ref{uni-sec}) of the solution of $(\sP)$
are characterised in terms of coverings and minimal 
coverings by sets which are inverses of subdifferentials of $g$.
As in the work of Zimmermann, we obtain an algorithm
to check the existence and uniqueness of the solution of $(\sP)$, 
when $X$ and $Y$ are finite (Section~\ref{gal-ex-sec}).
We also illustrate our results on various Moreau conjugacies
(Section~\ref{sec4}).
When $B$ is the Legendre-Fenchel transform,
our results show (see Section~\ref{fenc-sec})
that essentially smooth convex functions have a
unique preimage by the Legendre-Fenchel transform, a fact which is
the essence of the classical G\"artner-Ellis theorem,
see e.g.~\cite[Theorem~2.3.6,(c)]{DEMBO} for a general presentation
(the use of the uniqueness of the preimage of the Legendre-Fenchel transform
was made explicit by O'Brien and Vervaat~\cite[Theorem~4.1 (c)]{VERV95},
Gulinsky~\cite[Theorems 4.7 and 5.3]{gulinski}
and Puhalskii~\cite[Lemmas 3.2 and 3.5]{PUHAL94-1}).
In Section~\ref{lip-sec}, 
we consider the Moreau conjugacy with
kernel $b(x,y)=-\omega(x-y)$, where $\omega$
is a nonnegative, continuous, and subadditive map.
Spaces of Lipschitz or H\"older continuous maps
arise as images of such conjugacies.
In Section~\ref{lsc-sec}, we consider the 
Moreau conjugacy with
kernel $b(x,y)=-x'' \|x'-y\|^p$, where $x=(x',x'')$,
$x''\geq 0$ and $p>0$. This conjugacy
was already studied by Dolecki and Kurcyusz~\cite{dolecki}.

Problem $(\sP)$ arises when looking for the rate function
in large deviations:
further applications of our results are given in~\cite{AGK2},
where we use our characterisations of the uniqueness in Problem  $(\sP)$ 
to give a new proof, as well as generalisations, of 
the G\"artner-Ellis theorem.
A second motivation arises
from optimal control: reconstructing
the initial condition of an Hamilton-Jacobi
equation from the final value is a special
case of Problem~$(\sP)$, which has been
studied, under convexity assumptions,
by Goebel and Rockafellar~\cite{goebel}.
Finally, Problem $(\sP)$ arises in the 
characterisation of dual solutions of the Monge-Kantorovitch 
mass transfer problem (see~\cite{rachev,villani} for general presentations).
Subdifferentials associated to Moreau
conjugacies are instrumental in this theory, as shown
by R\"uschendorf~\cite{frechet,ruschendorf95} 
(see also~\cite[Section~3.3]{rachev}). 

The results of the present paper were announced in~\cite{AGK1}.

{\em Acknowledgement.} We thank Ivan Singer and
an anonymous referee for their helpful comments on a
preliminary version of this paper. In particular, Ivan Singer
pointed out previous works on generalised subdifferentials,
and the anonymous referee pointed out the seminal 
work of Vorobyev.

\section{Representation of Functional Galois Connections}\label{sec-rep}
Many basic results of convex analysis are specialisations
of general properties of Galois connections in lattices,
that we next recall.
Galois connections between lattices of subsets
were introduced by Birkhoff, in the initial edition of~\cite{BIRK}.
Galois connections between general lattices
were introduced by Ore~\cite{ore}. (The word Galois correspondence
is sometimes used as a synonym of Galois connection.)
The proofs of the following results can be found
in~\cite[Chapter V, Section 8]{BIRK},~\cite{DUB},%
~\cite[chapter 1, Section 2]{BLY},%
~\cite[Section 4.4.2]{BCOQ}, 
and~\cite[Chapter 0, Section 3]{LATT}.

Let $(\sF,\leqf)$ and $(\sG,\leqg)$ be two partially 
ordered sets, and 
let $B:\sF\rightarrow\sG$ and $C:\sG\rightarrow\sF$.
We say that $B$ is \NEW{antitone}
if $f\leqf f' \implies Bf'\leqg  Bf$.
The pair $(B, C)$ is a {\em dual Galois connection} 
between $\sF$ and $\sG$ if
it satisfies one of the following equivalent conditions:
\begin{subequations}
\label{galois}
\begin{gather}
\label{galdef1}
I_\sF \geq_\sF C B, \quad  I_\sG \geq_\sG B C ,
\quad \text{and}\; B,C\;\mrm{are antitone maps},\\
 \label{galdef2}
(g\geq_\sG B f \iff f\geq_\sF C g) \qquad
\forall f\in \sF, \;g\in \sG\enspace,\\
\label{resid}
C g=\minf \set{ f}{ g \geq_\sG  B f}\qquad \forall g\in \sG\enspace,\\
B f=\ming \set{ g}{f \geq_\sF  C g}  \qquad \forall f\in \sF\enspace,
\end{gather}
\end{subequations}
where $\minf$ (resp.\ $\ming$) denotes the minimum element of a set
for the order $\leqf$ (resp.\ $\leqg$),
and where $I_\AC$ denotes the identity on a set $\AC$.

It follows from~\eqref{resid}
that for any $B$, there is at most one map $C$ such that $(B,C)$ is 
a dual Galois connection. We denote this $C$ by $B\gal$.
It also follows from~\eqref{resid} that
for all $g\in \sG$ and $h\in \sF$,
\begin{equation}
g=Bh\impl B\gal g=\minf \set{ f}{g= B f }\enspace.
\label{ef}
\end{equation}
In particular, 
\begin{equation}
B\gal =B^{-1} \quad \text{\rm if $B$ is invertible.}
\label{e-bisinvertible}
\end{equation}
The two inequalities in~\eqref{galdef1} imply that
\begin{align*}
B B\gal B= B\qquad  \text{\rm and}\qquad B\gal B B\gal=B\gal\enspace.
\end{align*}
{From} this, or from~\eqref{ef}, one deduces
\begin{align*}
Bf=g \text{ has a solution } f\in \sF \quad \Longleftrightarrow 
\quad B B\gal g=g\enspace.
\end{align*}  
If $(B,C)$ is a dual Galois connection, so does $(C,B)$, by symmetry.
Hence, $(B\gal )\gal=B$. If $B$ yields a dual Galois connection,
then
\begin{equation}\label{infi}
B \left(\inff F\right) = \supg\set{Bf}{f\in F}\enspace,\end{equation}
for any subset $F$ of $\sF$ such that the infimum of $F$ exists.
In particular, if $\sF$ has a maximum element, $\topsf$,
we get by specialising~\eqref{infi} to $F=\emptyset$
that $\sG$ has a minimum element, $\bot_{\sG}$,
and
\begin{equation}
B(\topsf)= \bot_{\sG}\enspace.
\label{e-toptobot}
\end{equation}
Moreover, if $\sF$ is a complete ordered set, i.e.,
if any subset of $\sF$ has a greatest lower bound,
property~\eqref{infi} characterises
the maps $B$ that yield a dual Galois connection.

Ordinary (non dual) {\em Galois connections} are defined by 
reversing the order relation of $\sF$ and $\sG$ in~\eqref{galois}.
One also finds in the literature the names
of {\em residuated maps} $B$ and {\em dually residuated maps} $C$, which are
defined by reversing the order of $\sF$, but not the order
of $\sG$, in~\eqref{galois} (see for instance~\cite{BLY,BCOQ}).
All these notions are equivalent.

We call \NEW{lattice of functions}
a sublattice $\sF$ of $S^Y$, where
$(S,\leq)$ is a lattice, $Y$ is a set,
and $S^Y$ is equipped with the product ordering
(that we still denote by $\leq$).
When $\sF\subset S^Y$ and
$\sG\subset T^X$ are
lattices of functions,
we say that $(B,B\gal)$ is a (dual) \NEW{functional} 
Galois connection.  

When $S$ has a maximum element $\topS$, $y\in Y$ and $s\in S$,
we denote by $\delta^s_y$ the map:
\begin{align*}
\delta^s_y\in S^Y:\qquad 
\delta^s_y(y')=\begin{cases}
 s & \mrm{if}\; y'=y,\\
\topS& \mrm{otherwise,}
\end{cases}
\end{align*}
that we call the \NEW{Dirac function} at point $y\in Y$ with value $s\in S$.

\begin{theorem}\label{theo-rep}
Let $S,T$ be two lattices that have a maximum element,
let $X,Y$ be arbitrary nonempty sets
and let $\sF\subset S^Y$ (resp.\ $\sG\subset T^X$)
be a lattice of functions containing
all the Dirac functions of $S^Y$ (resp.\ $T^X$).
Then $(B,B\gal)$ is a dual Galois connection
between $\sF$ and $\sG$ if, and only if, 
there exist two maps $b: X\times Y\times S\to T$
and $\bgal: Y\times X \times T\to S$ such that:
for all $(x,y)\in X\times Y$, $(b(x,y,\cdot),\bgal(y,x,\cdot))$
is a dual Galois connection between $S$ and $T$;
for all $(x,t)\in X\times T$, $\bgal(\cdot,x,t)\in \sF$;
for all $(y,s)\in Y\times S$, $b(\cdot,y,s)\in \sG$;
and
\begin{subequations}
\label{e-repBB}
\begin{align}
Bf&= \supg \set{b(\cdot ,y,f(y))}{y\in Y}
,\quad \forall f\in \sF\enspace, \label{e-repB}\\
B\gal g&= \supf \set{\bgal(\cdot,x,g(x))}{x\in X} 
,\quad \forall g\in \sG\enspace.\label{e-repBgal}
\end{align}
\end{subequations}
In this case, the maps $b$ and $b\gal$ are uniquely determined by 
$(B,B\gal)$, since $b(\cdot,y,s)= B \delta_y^s$ and
$b\gal(\cdot,x,t)=B\gal \delta_x^t$ for all $s\in S, t\in T, x\in X$
and $y\in Y$.
\end{theorem}
Theorem~\ref{theo-rep}
was inspired by a ``Riesz representation theorem'' 
of Maslov and Kolokoltsov~\cite{KOLO88,KOLO90,kolokoltsov92}
(see also~\cite[Theorem~1.4]{maslovkolokoltsov95}) 
which is similar to Theorem~\ref{theo-rep}: it 
applies to a continuous map $B$ between (non-complete)
lattices of {\em continuous} functions
$\sF$ and $\sG$, with  $S=T=\rbar$,
assuming that $B$ preserves finite sups.
Mart\'\i nez-Legaz and Singer obtained 
in~\cite[Theorems~3.1 and 3.5]{martinez90} (see 
also~\cite[Theorem~7.3]{singer}) the same conclusions as in 
Theorem~\ref{theo-rep}, in the special case where $\sF=S^Y$, $\sG=T^X$ and
$S$ and $T$ are complete lattices 
(Theorem~7.3 of~\cite{singer} is stated in the 
case where $S$ and $T$ are included in $\rbar$, but
it is remarked in~\cite[page 419]{singer}  that this result 
is valid for general complete lattices $S$ and $T$).

Theorem~\ref{theo-rep} allows us to consider
the case where $\sF=\lsc(Y,S)$ is the set
of \NEW{lower semicontinuous}, or \NEW{l.s.c.}, maps from $Y$ 
to $S$. Here, we say that a map $f: Y\to S$ is
l.s.c.\ if for all $s\in S$, the sublevel set 
$\set{y\in Y}{f(y)\leq s}$ is closed.
When $Y$ is a $T_1$ 
topological space and $S$ has a maximum element,
the Dirac functions are l.s.c.,
so that Theorem~\ref{theo-rep} can be applied.
In this case, $\supf=\sup$ since the sup of l.s.c.\ maps is l.s.c.,
and Theorem~\ref{theo-rep} shows that $b\gal(\cdot,x,t)$
is l.s.c. We shall see in Proposition~\ref{prop-blsc} below
that $b(x,\cdot,s)$ is also l.s.c.
\begin{remark}
If $S,T,\sF,\sG$ are as in Theorem~\ref{theo-rep},
$\sF$ has a maximum element, namely, the constant function
$y\mapsto \topS$ (which necessarily belongs to $\sF$
because it is equal to the Dirac function $\delta_y^{\topS}$
for any $y\in Y$).
Then, if $(B,B\gal)$ is a dual Galois connection between 
$\sF$ and $\sG$, the remark before~\eqref{e-toptobot} shows that $\sG$
has a minimum element.
Moreover, by Theorem~\ref{theo-rep}, the existence of a dual
Galois connection between $\sF$ and $\sG$ implies
the existence of dual Galois connection between $S$ and $T$, 
hence, by \eqref{e-toptobot},
 $T$ has a minimum element,  $\bot_{T}$.
Since $\sG$ contains all the Dirac functions, the minimum element of $\sG$,
$\bot_\sG$, is such that $\bot_\sG(x)\leq \delta_x^{\bot_T}(x)=
\bot_T$ for all $x\in X$, hence $\bot_\sG$ 
is necessarily the constant function $x\mapsto \bot_T$.
Symmetrically, $S$ (resp.\ $\sF$)  has a minimum element,
$\bot_S$ (resp.\ the constant function $y\mapsto \bot_S$).
\end{remark}

\begin{proof}[Proof of Theorem~\ref{theo-rep}]
The proof of Theorem~\ref{theo-rep} is similar to that of
Theorems~3.1 and 3.5 in~\cite{martinez90}. We give it for completeness.
Let us first assume that $(B,B\gal)$
is a dual Galois connection, and define
\[
b(\cdot,y,s) = B \delta^s_y\in \sG,\qquad 
\bgal(\cdot,x,t) = B\gal \delta^t_x\in \sF\enspace.
\]
We have
\begin{align*}
b(x,y,s)\leq t & \iff B \delta^s_y (x)\leq t\\
&\iff B \delta^s_y\leq \delta^t_x\\
& \iff B\gal \delta^t_x \leq \delta^s_y \quad \mrm{(by~\eqref{galdef2})} \\
& \iff \bgal(y,x,t)\leq s\enspace, 
\end{align*}
which shows, by~\eqref{galdef2},
that $(b(x,y,\cdot),\bgal(y,x,\cdot))$
is a dual Galois connection between $S$ and $T$.

Using~\eqref{infi} and
$f=\inff\set{\delta_y^{f(y)}}{y\in Y}$,
which holds for all $f\in \sF$, we get~\eqref{e-repB}.
The representation~\eqref{e-repBgal}
is obtained by symmetry.

Conversely, let us assume that $(B,B\gal)$
are defined by~\eqref{e-repBB} 
where $b$ and $\bgal$ satisfy the conditions of the theorem.
(This means in particular
that the $\supg$ and $\supf$ in~\eqref{e-repBB} exist.)
Then, applying~\eqref{e-repB} to $\delta_y^s\in\sF$ and
using~\eqref{e-toptobot}, we get that $b(\cdot,y,s)= B \delta_y^s$.
Similarly, $b\gal(\cdot,x,t)=B\gal \delta_x^t$.
Moreover, for all $f\in \sF,g\in \sG$,
\begin{align*}
Bf \leq g &\iff
b(\cdot,y,f(y))\leq g,\; \forall y\in   Y\\
&\iff b(x,y,f(y))\leq g(x),\; \forall (x,y)\in  X\times Y\\
& \hspace{8em} \mrm{(since $\leq$ is the product ordering on $\sG$)}\\
&\iff \bgal(y,x,g(x))\leq f(y),\; \forall (x,y)\in  X\times Y
\;\;\mrm{(by~\eqref{galdef2})}\\
&\iff B\gal g \leq f\enspace,
 \end{align*}
which, by~\eqref{galdef2} again, shows
that $(B,B\gal)$ is a dual Galois connection.
\end{proof}
\begin{proposition}\label{prop-blsc}
If $(b,\bgal)$ are as in Theorem~\ref{theo-rep},
where $Y$ is a $T_1$ topological space
and $\sF=\lsc(Y,S)$, then, for all $(x,s)\in X\times S$, the map
$b(x,\cdot,s): Y\to T$ is l.s.c.
\end{proposition}
\begin{proof}
By Theorem~\ref{theo-rep},
$\bgal(\cdot,x,t)\in \sF=\lsc(Y,S)$,
for all $(x,t)\in X\times T$.
The equivalence~\eqref{galdef2}
shows that $\bgal(\cdot,x,t)$ 
is l.s.c.\ for all $t\in T$,  
if, and only if, 
$b(x,\cdot,s)$ is l.s.c.\ for all $s\in S$.
\end{proof}
By symmetry, when $X$ is a $T_1$ topological space,
and $\sG=\lsc(X,T)$, the map $\bgal(y,\cdot,t)$ is l.s.c.\ 
for all $(y,t)\in Y\times T$.

We say that $\sF$ is a {\em lattice of subsets}
if there exists a set $Y$ such that
 $\sF\subset \pY$, the set of all subsets of $Y$,
 and $\sF$ is a lattice for the $\subset$ ordering.
Taking for $S$ the complete lattice of Booleans $(\{0,1\},\leq)$,
we can identify $\sF$ to a lattice of functions included in $S^Y$,
by using the lattice isomorphism: $F\in\pY\mapsto 1_F$, where
$1_F(y)=1$ if $y\in F$ and $1_F(y)=0$ otherwise.
Hence, specialising Theorem~\ref{theo-rep} to the case where $S$ and $T$
are equal to the lattice of Booleans, and considering non dual
Galois connections, we get: 
\begin{corollary}\label{galois-ens}
Let $X,Y$ be arbitrary nonempty sets
and let $\sF\subset \pY$ {\rm (}resp.\ $\sG\subset \pX${\rm )}
be a lattice of subsets containing
the empty set and all singletons of $Y$ (resp.\ $X$).
Then $(B,B\gal)$ is a Galois connection
between $\sF$ and $\sG$ if, and only if, 
there exists a set $\sB\subset X\times Y$ such that:
for all $x\in X$, $\sB_x=\set{y\in Y}{(x,y)\in\sB}\in \sF$;
for all $y\in Y$, $\sB^y=\set{x\in X}{(x,y)\in\sB}\in \sG$;
and
\begin{subequations}
\label{e-repgal}
\begin{align}
BF&= \infg \set{\sB^y}{y\in F}
,\quad \forall F\in \sF\enspace,\\
B\gal G&= \inff \set{\sB_x}{x\in G} 
,\quad \forall G\in \sG\enspace.
\end{align}
\end{subequations}
In this case, the set $\sB$ is uniquely determined by 
$(B,B\gal)$, since $\sB=\cup_{y\in Y} B\left(\{y\}\right)\times \{y\}
 =\cup_{x\in X} \{x\}\times B\gal \left(\{x\}\right)$.
\end{corollary}
When $\sF$ (resp.\ $\sG$) is stable by taking arbitrary intersections,
$\inff$ (resp.\ $\infg$) coincides with the intersection operation.

The conclusions of Corollary~\ref{galois-ens} were obtained
by Everett, 
when $\sF\subset \pY$ and $\sG\subset\pX$ 
are complete distributive lattices
(see Theorem~5 of~\cite{everett},
the remark following its proof,
and Section 8 of Chapter 5, page 124 of~\cite{BIRK}).
In the special case where $\sF=\pY$ and $\sG=\pX$,
the conclusions of Corollary~\ref{galois-ens}
were obtained in~\cite[Theorem 1.1]{singer86}
and~\cite[Theorem~3.3 and Remark 3.2]{martinez90}
using ideas of abstract convex analysis.

\begin{remark}
Many classical Galois connections are of the form~\eqref{e-repgal}
(but $\sF$ and $\sG$ need not contain the singletons and 
the empty set).
For instance, if $Y$ is an extension of a field $K$,
if $\sF$ is the set of intermediate fields
$F$: $K\subset F\subset Y$, if $X$ is the group of automorphisms of $Y$ fixing
every element of $K$, and if $\sG$ is the set of
subgroups of $X$, we obtain the original
Galois correspondence by setting $\sB=\set{(g,y)\in X\times Y}{g(y)=y}$.
\end{remark}

In the sequel, we shall only consider
the case where $S=T=\rbar$.
In this case, the property that $(b(x,y,\cdot),\bgal(y,x,\cdot))$
is a dual Galois connection can be made explicit:
\begin{lemma}\label{galoisr}
A map $h: \rbar \to \rbar$ yields a dual Galois
connection if, and only if, $h$ is nonincreasing,
right-continuous, and $h(+\infty)=-\infty$.  
\end{lemma}
\begin{proof}
This follows readily from the characterisation~\eqref{infi}
of dual Galois connections between complete lattices.
\end{proof}
Since $h$ is nonincreasing, one can replace right-continuous by l.s.c.\
in Lemma~\ref{galoisr} as is done in the statement 
of~\cite[Theorem 3.2]{martinez90}.

\begin{example}
When $b\in\rbar$, the map $h:\rbar\to\rbar,\; \lambda\mapsto b-\lambda$,
with the convention that $-\infty$ is absorbing
for addition, yields a dual Galois connection (by Lemma~\ref{galoisr}).
Moreover $h\gal=h$.
\end{example}
\begin{example}\label{exmoreau}
Let $S=T=\rbar$, and $X,Y,\sF,\sG$ be  as in Theorem~\ref{theo-rep}.
Assume in addition 
that $\sF$ and $\sG$ are stable by the addition of a constant,
again with the convention that $-\infty$ is absorbing for addition.
Let $\bar{b}: X\times Y\to \rbar$ be a map, and let $B:\sF\to \sG$ and
$B\gal:\sG\to\sF$ be defined by~\eqref{e-repBB} with
\begin{equation}\label{defib}
b(x,y,\alpha)=\bgal(y,x,\alpha)=\bar{b}(x,y)-\alpha
\quad\forall x\in X, \; y\in Y, \; \alpha\in\rbar\enspace.
\end{equation}
Theorem~\ref{theo-rep} shows that $(B,B\gal)$ is a dual Galois connection
if, and only if, $\bar b(x,\cdot)\in \sF$ for all $x\in X$,
and $\bar b(\cdot,y)\in \sG$ for all $y\in Y$.
This result can be applied, in particular, when $\sF=\lsc(Y,\rbar)$ and
$\sG=\RBX$. Such functional Galois connections are called
\NEW{Moreau conjugacies}~\cite{moreau70} and have been considered
by several authors (see for instance~\cite{dolecki} 
and~\cite[Section 5]{martinez90}).
In particular, taking two topological vector spaces in duality
$X$ and $Y$, and $\bar{b}: X\times Y\to \rbar,\; (x,y)\mapsto
\<x,y>$, we obtain the classical Legendre-Fenchel transform
$f\mapsto Bf=f^{\star}$.
\end{example}

\begin{remark}
The set $\rbar$ can be equipped with the semiring
structure of $\rmaxb$, in which the addition
is $(a,b)\mapsto \max(a,b)$ and the multiplication is
$(a,b)\mapsto a+b$, with the convention that $-\infty$
is absorbing for the multiplication of this semiring.
Then, if $Z$ is a set, $\RBZ$ can be equipped with two different
$\rmaxb$-semimodule structures. The \NEW{natural semimodule},
denoted $\RBBZ$, is obtained by taking the addition
$(f,f')\mapsto f\oplus f'$, with $(f\oplus f')(z)=\max(f(z),f'(z))$
for all $z\in Z$, and the action
$(\lambda,f)\mapsto \lambda. f$ with $(\lambda.f)(z)=\lambda +f(z)$
for all $z\in Z$, again with the convention that $-\infty$
is absorbing. 
The \NEW{opposite semimodule}, denoted  $(\RBBZ)^{\mrm{op}}$, 
is obtained by taking the addition
$(f,f')\mapsto f\oplus' f'$, with $(f\oplus' f')(z)=\min(f(z),f'(z))$
for all $z\in Z$, and the action
$(\lambda,f)\mapsto \lambda.' f$ with $(\lambda.' f)(z)=-\lambda +f(z)$
for all $z\in Z$, with the dual convention
that $+\infty$ is absorbing (see~\cite{cgq02}).
Then the Moreau conjugacies, i.e., the functional Galois  connections
of Example~\ref{exmoreau}, are $\rmaxb$-linear from $(\RBBZ)^{\mrm{op}}$
to $\RBBZ$.
\end{remark}

\section{Existence of Solutions of $Bf=g$}\label{exis-sec}

\subsection{Statement of the Existence Result}\label{sec31}
In the following, we take $S=T=\rbar$, 
we assume that $X$ and $Y$ are
Hausdorff topological spaces, and take
$\sF=\lsc(Y,\rbar)$, $\sG=\RBX$,
together with $B,B\gal,b,\bgal$ as in Theorem~\ref{theo-rep}.
This includes the case where $\sF=\RBY$ and $X,Y$ are arbitrary sets,
which is obtained by taking discrete topologies on $X$ and $Y$.
Since, by Theorem~\ref{theo-rep}, $(b(x,y,\cdot),\bgal(y,x,\cdot))$
is a dual Galois connection, Lemma~\ref{galoisr} shows that
$b(x,y,\cdot)$ and $\bgal(y,x,\cdot)$ are nonincreasing and right-continuous
maps from $\rbar$ to itself, which take the value $-\infty$ at $+\infty$.

We shall assume in the sequel that there is
a subset $\sS\subset X\times Y$ satisfying:
\begin{itemize}
\newcounter{assume}
\def\theassume{A\arabic{assume}}
\def\myitem{\refstepcounter{assume}\item[(\theassume)]}
\myitem\label{a1}
$\sS_x=\set{y\in Y}{(x,y)\in \sS}\neq\emptyset$,
for all $x\in X$; 
\myitem\label{as2}
$\sS^y=\set{x\in X}{(x,y)\in \sS}\neq\emptyset$,
for all $y\in Y$; 
\myitem\label{a3} $b(x,y,\cdot)$ is a bijection
$\rbar\to \rbar$ for all $(x,y)\in \sS$;
\myitem\label{assump4} $b(x,y,\cdot)\equiv -\infty$, for $(x,y)\in
X\times Y\setminus \sS$.
\end{itemize}

When $B$ is the Moreau conjugacy given by~{\rm (\ref{e-repBB},\ref{defib})},
Assumptions~{\rm (}\ref{a1}--\ref{assump4}{\rm )}
 are satisfied if, and only if,
 $\bar{b}(x,y)\in\R\cup\{-\infty\}$ for all $(x,y)\in X\times Y$, 
and $\sS:=\set{(x,y)\in X\times Y}{ \bar{b}(x,y)\in \R}$
satisfies~{\rm (}\ref{a1}--\ref{as2}{\rm )},
that is, for all $x\in X$ and $y\in Y$, $\bar{b}(x,\cdot)$
and $\bar{b}(\cdot,y)$ are not identically $-\infty$.
These assumptions are fulfilled in particular for the kernel
of the Legendre-Fenchel transform.

Rather than Problem $(\sP)$, we will
consider the more general problem:
\[
(\sP'):\qquad \mrm{Find }f\in\sF\mrm{ such that } Bf\leq g
\mrm{ and }  Bf(x)=g(x)\mrm{ for all } x\in X'\enspace,
\]
where $g\in \sG$ and $X'\subset X$ are given.

To state our results, we need some definitions and notations.
When $B$ is the Legendre-Fenchel transform, these notations and definitions 
will correspond to those defined classically in convex analysis.
First, for any map $g$ from a topological space $Z$ to $\rbar$,
we define the \NEW{lower domain}, \NEW{upper domain},
\NEW{domain}, and \NEW{inner domain}:
\begin{align*}
\ldom g &= \set{z\in Z}{g(z)<+\infty}\enspace, \\
\udom g &= \set{z\in Z}{g(z)>-\infty}\enspace, \\
\dom g  &= \ldom g \cap \udom g\enspace,  \\
\idom g &=  \set{z\in \dom g}{\displaystyle
\limsup_{z'\to z} g(z')<+\infty}\enspace.
\end{align*}
The set $\idom g\subset \dom g$ is an open subset of $\udom g$. 
When $Z$ is endowed with the discrete topology, $\idom{g}=\dom{g}$. 
We say that $g$ is \NEW{proper} if $g(z)\neq -\infty$ for all $z\in Z$ and
if there exists $z\in Z$ such that $g(z)\neq +\infty$, which means that
$\udom{g}=Z$ and $\dom g\neq\emptyset$.

We shall use the following variant of the 
general notion of subdifferentials of dualities
introduced by Mart\'\i nez-Legaz and Singer~\cite{martinez95}
(see Remark~\ref{rem-sing} for a comparison). 
Given  $f\in\sF$ and $y\in Y$, we call 
the \NEW{subdifferential} of $f$ at $y$ with respect to $b$ (or $B$), 
and we denote by $\partial^b f(y)$, or $\partial f(y)$ for brevity,
the set:
\begin{subequations}
\label{def-subgen}
\begin{align}\label{defpart}
\partial f(y)&= \set{x\in X}{(x,y)\in \sS,\quad
b(x,y',f(y')) \leq b(x,y,f(y))\; \forall y'\in Y}\enspace.
\end{align}
For $g\in \sG$, and $x\in X$, the subdifferential of $g$ at $x$
with respect to $b\gal$,  $\partial^{b\gal} g(x)$, that 
we denote $\partial\gal g (x)$ for brevity, is given by:
\begin{align}
\label{defpartr}
\partial\gal g(x)&= \set{y\in Y}{(x,y)\in\sS,\quad
\bgal(y,x', g(x'))\leq \bgal(y,x,g(x))\; \forall x'\in X}\enspace.
\end{align}
\end{subequations}
Then 
\begin{subequations}
\label{e-inverse}
\begin{align}
\partial f(y)&= \set{x\in X}{(x,y)\in\sS\text{ and }
Bf(x)= b(x,y,f(y))}\enspace,\label{e-inverse1}\\
\partial\gal g(x)&= \set{y\in Y}{(x,y)
\in\sS\text{ and } B\gal g(y)= \bgal(y,x,g(x))}\label{e-inverse2}\enspace,
\end{align}
\end{subequations}
and when $b(x,y,\alpha)=\<x,y>-\alpha$ is the kernel
of the Legendre-Fenchel transform,
we recover the classical definition
of subdifferentials. 
\begin{remark}\label{rem-sing}
{From}~\eqref{a3} and~\eqref{e-bisinvertible}, 
$b(x,y,\cdot)$ is a bijection with inverse $\bgal(y,x,\cdot)$,
for all $(x,y)\in\sS$.  
Hence,
\begin{align*}
\partial f(y)&= \set{x\in X}{(x,y)\in\sS\text{ and }
b\gal (y,x,Bf(x))= f(y)}\enspace,\\
\partial\gal g(x)&= \set{y\in Y}{(x,y)
\in\sS\text{ and } b(x,y,B\gal g(y))=g(x)}\enspace.
\end{align*}
In ~\cite{martinez95}, the subdifferential $\partial^B f(y)$
of $f$ at $y$ with respect to $B$ is defined by $\partial^B f(y)
= \set{x\in X}{b\gal (y,x,Bf(x))= f(y)}$. Hence,
under  Assumptions~{\rm (\ref{a1}--\ref{assump4})},
the subdifferential $\partial f(y)$ and the subdifferential 
$\partial^B f(y)$ of Mart\'\i nez-Legaz and Singer 
differ only when $f(y)=-\infty$.
In this case $\partial^B f(y)= X$ and $\partial f(y)=\sS^y$.
Hence, $\partial f(y)\subset \partial^B f(y)$
holds for all $y\in Y$.
Moreover $\partial^B f(y)=\emptyset$ if, and only if, 
$\partial f(y)=\emptyset$. With these remarks
in mind, one may rephrase equivalently all the subsequent results
in terms of $\partial^B$ and $\partial^{B\gal}$. 
\end{remark}

As was observed in~\cite{martinez95}, the
generalised subdifferentials have a geometric
interpretation.
Let us fix $(x,y)\in\sS$.
Thanks to Assumption~\eqref{a3},
there is a unique function of the family 
$\{b(\cdot,y,\lambda)\}_{\lambda\in\rbar}$ which takes
the value $g(x)$ at point $x$: 
this function is obtained for $\lambda=\bgal(y,x,g(x))$,
since, as observed in Remark~\ref{rem-sing},
$\bgal(y,x,\cdot)$ is the inverse of $b(x,y,\cdot)$.
Let us call this function the {\em curve of direction $y$
meeting $g$ at point $x$}. 
By~\eqref{galdef2}, $y\in\partial\gal g(x)$
if, and only if, $g(x')\geq b(x',y,\bgal(y,x,g(x)))$,
for all $x'\in X$,
which means that {\em the curve of direction
$y$ meeting $g$ at $x$ is below $g$}.
When $B=B\gal$ is the Legendre-Fenchel transform,
the curves become lines, and we recover
the classical interpretation of subdifferentials.

\begin{example}\label{ex-geom}
The geometrical interpretation of subdifferentials is illustrated
in Figure~\ref{fig2}, where $X=Y=\R$, $b$ is given by~\eqref{defib}
with $\bar{b}(x,y)=-|x-y|$,
$g(x)=x^2/2$ for $x\leq 0$,
$g(x)=x$ for $x\in [0,1]$, $g(x)=1$ for $x=[1,3)$, 
and $g(x)=x/3-1$ for $x\in[3,\infty)$. We have $\partial\gal g(x)=\emptyset$
for $x\in (-\infty,-1)\cup (2,3)$, 
$\partial \gal g(x)=\{x\}$ for $x\in (-1,0)\cup [1,2]\cup (3,\infty[$,
$\partial\gal g(x)=[x,1]$ for $x\in [0,1)$,
$\partial\gal g(-1)=(-\infty,-1]$ and $\partial\gal g(3)=[2,3]$.
This can be checked
by looking at the curves meeting $g$ at points
$-1$, $1/2$, $3$ and $4$, which are depicted on the figure.
\begin{figure}[htb]
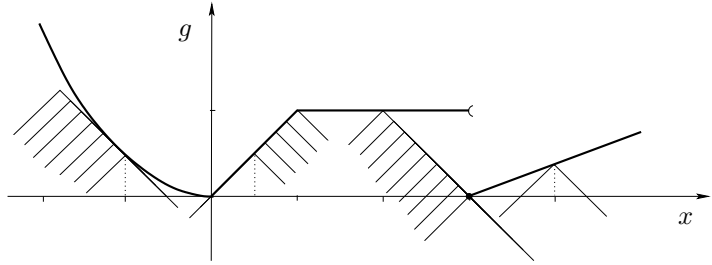

\input fig2
\caption{Geometric interpretation of subdifferentials}
\label{fig2}
\end{figure}
\end{example}
\begin{definition}\label{defi-cover}
When $F$ is a map from a set $Z$ to the set $\pW$ of all subsets of some
set $W$, we denote by $F^{-1}$ the map from $W$ to $\pZ$ given by
$F^{-1}(w)=\set{z\in Z}{w\in F(z)}$,
and we define the domain of $F$: 
$\dom{F}:=\set{z\in Z}{F(z)\neq\emptyset}=\cup_{w\in W} F^{-1}(w)$.
If $Z'\subset Z$ and $W'\subset W$,
we say that $\{F(z)\}_{z \in Z'}$ is a {\em covering}
of $W'$ if $\cup_{z\in Z'} F(z)\supset W'$. 
\end{definition}
When $F,Z,Z',W$ are as in Definition~\ref{defi-cover},
the family $\{F^{-1}(w)\}_{w \in W}$ is a covering
of $Z'$ if, and only if, $Z'\subset \dom{F}$.
By~\eqref{def-subgen},
\begin{subequations} 
\begin{align}
(\partial f)^{-1}(x)&= \argmax_{y\in \sS_x} b(x,y,f(y))\enspace,
\\
(\partial\gal g)^{-1}(y)&= \argmax_{x\in \sS^y} \bgal(y,x,g(x))
\label{argmax2}\enspace.
\end{align}
\end{subequations}

\begin{definition}
We say that $b$ is \NEW{continuous in the second variable} if for all $x\in X$
and $\alpha\in \R$,
$b(x,\cdot,\alpha)$ is continuous.
We say that $b$ is \NEW{coercive} 
if for all $x\in X$, all neighbourhoods $V$ of $x$ in $X$,
and all $\alpha\in \R$, the function 
\begin{equation}
y\in Y\mapsto b_{x,V}^{\alpha}(y)=
\sup_{z\in V} b(z,y,\bgal(y,x,\alpha))\enspace, 
\label{e-def-bv}
\end{equation}
has relatively compact finite sublevel sets, which
means that $\set{y\in Y}{b_{x,V}^{\alpha}(y)\leq \beta}$
is relatively compact for all $\beta\in \R$. 
\end{definition}
The continuity of $b$ in the second variable
holds readily when $Y$ is discrete (and in particular when $Y$ is finite).
The coercivity of $b$ holds trivially,
and independently of the topology on $X$,
when $Y$ is compact (and in particular when $Y$ is finite).
If $X=Y=\R^n$ and $b(x,y,\alpha)=\<x,y>-\alpha$
then $b$ is continuous in the second variable, and 
for all  neighbourhoods $V$ of $x$,
and all $\alpha\in \R$,
$b_{x,V}^{\alpha}(y)\geq \varepsilon \|y\|+\alpha $,
for some $\varepsilon>0$, so that $b$ is coercive.
Similarly, if $b(x,y,\alpha)=a\|x- y\|^2-\alpha$,
where $a\in \R\setminus \{0\}$ and $\|\cdot\|$ is the Euclidean norm,  
then $b$ is continuous in the second variable, and 
for all  neighbourhoods $V$ of $x$,
and all $\alpha\in \R$,
$b_{x,V}^{\alpha}(y)\geq \varepsilon \|y-x\| -1 +\alpha $,
for some $\varepsilon>0$, so that $b$ is coercive.

We also denote by $\sFC$ the set of all $f\in \sF$ such that
for all $x\in X$, $y\mapsto b(x,y,f(y))$ has
relatively compact finite superlevel sets, which
means that for all $\beta\in\R$, the set 
$\set{y\in Y}{b(x,y,f(y))\geq \beta}$ is relatively compact.
When $Y$ is compact, $\sFC=\sF$.

We shall occasionally make the following assumptions:
\begin{itemize}
\catcode`\@=11
\def\refcounter#1{\protected@edef\@currentlabel
       {\csname the#1\endcsname}%
}
\catcode`\@=12
\def\theassume{{\rm (A\arabic{assume})}}
\def\theassumep{{\rm (A\arabic{assume})$'$}}
\def\myitem{\refstepcounter{assume}\item[\theassume\hskip 0.72ex]}
\def\myitemp{\refcounter{assumep}\item[\theassumep]}
\myitem\label{a-yisdiscrete} $Y$ is discrete;
\myitemp\label{a-bgalgisfinite} $b$ is continuous in the second variable,
and $B\gal g(y)>-\infty$ for all $y\in Y$; 
\myitem\label{a-bgalginfc}\label{casei} $B\gal g\in\sFC$;
\myitemp\label{caseii}
\label{a-biscoercive}
$b$ is coercive and $X'\subset \idom g\cup g^{-1}(-\infty)$.
\end{itemize}
The assumption that  $B\gal g(y)>-\infty$ for all $y\in Y$
is fulfilled in very general situations: if  $g(x)<+\infty$ for all $x\in X$,
in particular if $\idom g=\udom g$ as in Corollary~\ref{cor-exis2}
below; or if $\sS=X\times Y$
and  $g\not\equiv +\infty$, which is the case
for instance when $B$ is the Legendre-Fenchel transform
and $g$ is proper.

The following general existence 
result is proved in Section~\ref{sec-proofcover}.
\begin{theorem}\label{cover}
Let $X'\subset X$, and $g\in \sG$.
Consider the following statements:
\begin{align}
&\mrm{Problem }(\sP')\mrm{ has a solution,}
\label{e-eq1} \\
& X'\subset \dom{\partial\gal g}, \label{e-c0}\\
&\{(\partial\gal g)^{-1}(y)\}_{y\in Y} \mrm{ is a covering of } X',
\label{e-c1}\\
&\{(\partial\gal g)^{-1}(y)\}_{y\in \ldom{B\gal g}}
\mrm{ is a covering of } X' \cap \udom g, \label{e-c2}\\
&\{(\partial\gal g)^{-1}(y)\}_{y\in \dom{B\gal g}}
\mrm{ is a covering of } X' \cap \dom g \label{e-c3}. 
\end{align}
We have: \eqref{e-c0}$\Leftrightarrow$\eqref{e-c1}$\Leftrightarrow$%
\eqref{e-c2}$\Rightarrow${\rm (\ref{e-eq1},\ref{e-c3})}.
The implication \eqref{e-eq1}$\Rightarrow$\eqref{e-c2} holds
if Assumptions~
\ref{a-yisdiscrete} or~\ref{a-bgalgisfinite}, 
and \ref{a-bgalginfc} or~\ref{a-biscoercive}, hold.
This is the case in particular if $Y$ is finite.
The implication \eqref{e-c3}$\Rightarrow$\eqref{e-c2} holds 
when $X'\subset \idom g\cup g^{-1}(-\infty)$.
Finally,~{\rm (\ref{e-eq1}--\ref{e-c3})} are true 
when $g\equiv +\infty$ or $g\equiv -\infty$.
\end{theorem}
The most intuitive condition should be~\eqref{e-c0}:
it says that $\partial\gal g(x)\neq\emptyset$,
for all $x\in X'$. See Example~\ref{ex-new} below
for an illustration. We stated conditions involving
coverings to make it clear that Theorem~\ref{cover} generalises 
the theorem of Vorobyev and Zimmermann. Moreover, coverings will
be instrumental in the statement of the uniqueness
results in Section~\ref{sec-unique}. 

\begin{remark}
As we shall see in Section~\ref{sec-proofcover}, the implication
\eqref{e-c0}$\Rightarrow$\eqref{e-eq1} can be deduced from~\eqref{tp3}.
The implication in~\eqref{tp3} was already shown by 
Mart\'\i nez-Legaz and Singer~\cite[Proposition 1.2]{martinez95},
using the notion of subdifferential of~\cite{martinez95}.
\end{remark}
We next state some direct corollaries.
\begin{corollary}\label{cor-exis3}
Consider $g\in \sG$. Assume that $Y$ is finite.
Then $Bf=g$ has a solution $f\in  \sF$
if, and only if, $\{( \partial\gal g)^{-1}(y)\}_{y\in Y}$ 
is a covering of $X$. \qed
\end{corollary}

\begin{corollary}\label{cor-exis2}
Consider $g\in \sG$ such that $\idom g=\udom g$.
Assume that $b$ is continuous in the second variable and coercive.
Then $Bf=g$ has a solution $f\in  \sF$
if, and only if, $\{( \partial\gal g)^{-1}(y)\}_{y\in Y}$ 
is a covering of $X$.\qed
\end{corollary}
\begin{example}
When $B$ is the Legendre-Fenchel transform over $\R^n$,
and $g$ is a l.s.c.\ proper convex function, 
Problem $(\sP)$ has a solution. A fortiori,
Problem $(\sP)'$ has a solution with $X'=\idom g$.
Then the implication \eqref{e-eq1}$\Rightarrow$\eqref{e-c0}
of Theorem~\ref{cover} shows that
$g$ admits subdifferentials in $\idom g$,
a well known result since for any l.s.c.\ convex function $g$ on $\R^n$,
$\idom g$ is the interior of $\dom g$ (see for 
instance~\cite[Theorem 23.4]{ROCK}).
Corollary~\ref{cor-exis2} shows that
if $g$ is finite and locally bounded from above everywhere,
then $g$ is l.s.c.\ and convex if, and only if,
it has nonempty subdifferentials everywhere.
\end{example}
The following final corollary is proved in Section~\ref{sec-proofcover}:
\begin{corollary}\label{cor-compact}
Let  $g\in \sG$.
Make Assumptions~
\ref{a-bgalgisfinite}, and \ref{a-bgalginfc}.
Then $Bf=g$ has a solution $f\in \sF$ if, and only if, 
$\{( \partial\gal g)^{-1}(y)\}_{y\in Y}$
is a covering of $X$. In that case, $g(x)<+\infty$
for all $x\in X$.
\end{corollary}
\begin{remark}\label{rem-a7tierce}
Assumptions~\ref{a-yisdiscrete} or~\ref{a-bgalgisfinite}, 
together with~\ref{a-bgalginfc} or~\ref{a-biscoercive} imply that:
\begin{itemize}
\catcode`\@=11
\def\refcounter#1{\protected@edef\@currentlabel
       {\csname the#1\endcsname}%
}
\catcode`\@=12
\def\theassume{(A\arabic{assume})}
\def\myitem{\refstepcounter{assume}\item[\theassume\hskip 0.72ex]}
\myitem\label{a7-tierce}
$\forall x\in X'\cap \udom g$, $\exists y\in Y$,
$b(x,y,B\gal g(y))=\sup_{y'\in Y}b(x,y',B\gal g(y'))$.
\end{itemize}
Indeed, Assumption~\ref{a-yisdiscrete} or~\ref{a-bgalgisfinite}
implies that the map $y\mapsto b(x,y,B\gal g(y))$ is u.s.c.\ 
for all $x\in X$, whereas  Assumption~\ref{a-bgalginfc} (resp.\
\ref{a-biscoercive}) implies that the map
$y\mapsto b(x,y,B\gal g(y))$ has relatively compact finite superlevel sets,
for all $x\in X$ (resp.\ for all $x\in \idom g\supset X'\cap \udom g$).
One can check that the conclusions of Theorem~\ref{cover} remain
valid if Assumptions~\ref{a-yisdiscrete}
or~\ref{a-bgalgisfinite}, and~\ref{a-bgalginfc} or~\ref{a-biscoercive},
are replaced by the single Assumption~\ref{a7-tierce}. 
Indeed, the proof of Theorem~\ref{cover} is based on the
application of Theorem~\ref{exis} to $f=B\gal g$, and
on the inclusion of $X'\cap\udom g$ in the set $X_0$ which 
appears in~\eqref{e-toprove}, and the proof of Theorem~\ref{exis} 
is precisely based on~\eqref{e-toprove}
(the other arguments do not need any assumption).
Similarly, the proof of Corollary~\ref{cor-compact} remains valid
if all the assumptions are replaced by Assumption~\ref{a7-tierce} 
with $X'\cap\udom g$ replaced by $X$.
\end{remark}
We shall give examples of application of Theorem~\ref{cover}
(and of its corollaries) in Sections~\ref{gal-ex-sec} and~\ref{sec4} below.
We next give examples illuminating the role of the technical
assumptions in Theorem~\ref{cover}.
\begin{example}\label{ex-new}
We first show a limitation of Theorem~\ref{cover}.
Consider again the kernel $b$ and the map 
$g$ of Example~\ref{ex-geom}. 
We claim that:
\begin{align}
\mrm{Problem }(\sP')\mrm{ has a solution, if, and only if, }
 X' \subset [-1,2]\cup[3,\infty)\enspace.
\label{e-cex}
\end{align}
Indeed, $ \dom{\partial\gal g}=\set{x\in\R}{\partial\gal g(x)\neq\emptyset}
=[-1,2]\cup[3,\infty)$, so that the ``if'' part of~\eqref{e-cex}
follows from the implication~\eqref{e-c0}%
$\Rightarrow$\eqref{e-eq1} in Theorem~\ref{cover}.
However, the ``only if'' part of~\eqref{e-cex} 
cannot be derived from Theorem~\ref{cover},
because Assumptions~\ref{a-bgalginfc} and~\ref{a-biscoercive} 
do not hold
(Assumption
~\ref{a-bgalgisfinite} is satisfied).
Indeed, we shall see in the more general setting
of Section~\ref{lip-sec} that the kernel $b$ is not coercive.
Moreover, one has  $B\gal g=- B B\gal g$, 
$B B\gal g(x)=g(x)$ for $x\in [-1,2]\cup[3,\infty)$,
$B B\gal g(x)=-x-1/2$ for $x\in (-\infty ,-1]$ and
$B B\gal g(x)=3-x$ for $x\in (2,3)$.
Then $b(x,y,B\gal g(y))=-x-1/2$ for $y\leq x$ and $y\leq -1$, and
$b(x,y,B\gal g(y))$ goes to $-\infty$ when $y$ goes to $+\infty$,
which implies that $B\gal g\not\in\sFC$.
Therefore, we cannot apply the implication 
\eqref{e-eq1}$\Rightarrow$\eqref{e-c2} in Theorem~\ref{cover}
to characterise the cases where Problem $(\sP')$ has a solution. 
However, if Problem $(\sP')$ has a solution, 
then by Lemma~\ref{lem-exis} below
and~\eqref{galdef1},
$B\gal g$ is also a solution of $(\sP')$. Therefore,
$X'$ must be included in the set
$\set{x\in X}{BB\gal g(x)=g(x)}$, which
is equal to $[-1,2]\cup[3,\infty)$,
by the previous computations.
This shows the ``only if'' part of~\eqref{e-cex}.
This conclusion can also be obtained by using the arguments of
Remark~\ref{rem-a7tierce} and the fact that 
Assumption~\ref{a7-tierce} holds for all subsets $X'$ of $\R$.
Note also that Theorem~\ref{cover} characterises
the solvability of Problem $(\sP')$ for any map $g\in\sG$
such that $\lim_{|x|\to\infty} g(x)-|x|=-\infty$,
since in this case, $B\gal g\in \sFC$.
\end{example}
\begin{example}
The following counter-example shows that the compactness
Assumptions~\ref{a-bgalginfc} or~\ref{a-biscoercive} are useful. 
Consider $X=\R$, $Y=[1,\infty)$, 
the Moreau conjugacy given by~{\rm (\ref{e-repBB},\ref{defib})}
with $\bar b(x,y)=0\vee (-|x|+ 1/y)$, where
$\vee$ denotes the sup law, and take the identically zero function $g$.
Then $B\gal g(y)=1/y$ and $g=B B\gal g$, 
which means that Problem  $(\sP')$ has a solution
with $X'=X$. However, $\partial\gal g(0)=Y$ and
$\partial\gal g(x)=\emptyset$ for all $x\in\R\setminus\{0\}$,
which shows that $\dom{\partial\gal g}=\{0\}$, hence
the covering condition~\eqref{e-c0} or~\eqref{e-c1} does not hold.
Since Assumption
~\ref{a-bgalgisfinite} is satisfied and
$\dom g=\idom g=X$, Theorem~\ref{cover} implies that $b$ is not coercive
and $B\gal g\not\in\sFC$.
\end{example}
\begin{example}\label{cex}
The following counter-example shows that Assumption~\eqref{a3}
is useful. Consider $Y=\{y_1\}$, $\sF=\RBY$,
$X=\{x_1,x_2\}$, $\sG=\RBX$, and the Moreau 
conjugacy given by~{\rm (\ref{e-repBB},\ref{defib})}
with $\bar{b}(x_1,y_1)=0$, $\bar{b}(x_2,y_1)=+\infty$.
We have $Bf(x_1)= -f(y_1)$, 
$Bf(x_2)=+\infty -f(y_1)$,
and $B\gal g(y_1)= \max( -g(x_1), (+\infty)-g(x_2))$,
for all $f\in\sF$ and $g\in\sG$.
Take $g(x_1) =0$ and 
$g(x_2) = +\infty $.
Then $Bf=g$ has a solution, namely, $f(y_1)=0$. 
However, taking any $\sS\subset X\times Y$ in~\eqref{argmax2}, we get
\(
(\partial\gal g)^{-1}(y_1)\subset
\argmax_{x\in \{x_1,x_2\}} \bar{b}(x,y_1)-g(x)= \{x_1\}
\)
which does not cover $X=\{x_1,x_2\}$. Therefore, 
the implication \eqref{e-eq1}$\Rightarrow$\eqref{e-c1} 
of Theorem~\ref{cover} does not extend to the case
of kernels $b$ which take the value $+\infty$, even when $Y$ is
finite (these kernels do not satisfy Assumption~\eqref{a3}).

Take now
$g(x_1) =-\infty  $
and $g(x_2) =0$.
Then $Bf=g$ has no solution, but taking $\sS= X\times Y$
in~\eqref{argmax2}, we get 
\(
(\partial\gal g)^{-1}(y_1)=
\argmax_{x\in \{x_1,x_2\}} \bgal(y_1,x,g(x))= X
\).
Therefore, the implication \eqref{e-c1}$\Rightarrow$\eqref{e-eq1}
of Theorem~\ref{cover} does not extend to the case of kernels
$b$ which take the value $+\infty$.
In this case, one should use rather the definition of subdifferentials
of Mart\'\i nez-Legaz and Singer, for which the
implication~\eqref{e-c1}$\Rightarrow$\eqref{e-eq1} can be deduced 
from~\cite[Proposition~1.2]{martinez95}.
\end{example}
\subsection{Additional Properties of $B$, and Proof of Theorem~\ref{cover}}
\label{sec-proofcover}

In this section, we state several lemmas and prove successively the
different assertions of Theorem~\ref{cover}.
We first show some properties of the kernels $b$ and $\bgal$.

By Theorem~\ref{theo-rep}, we know that
$(b(x,y,\cdot),\bgal(y,x,\cdot))$ is a {\em dual} Galois connection
between $\rbar$ and $\rbar$. The following result, which uses
Assumptions~\eqref{a3} and~\eqref{assump4}, shows
that $(b(x,y,\cdot),\bgal(y,x,\cdot))$ is almost a {\em (non dual)}
Galois connection:
\begin{lemma}\label{galoisnondual}
For all $(x,y)\in \sS$, $b(x,y,\cdot)$ is a decreasing
bijection $\rbar\to\rbar$ with inverse $\bgal(y,x,\cdot)$.
For all $(x,y)\in X\times Y$, we have
\begin{subequations}
\begin{equation}
\label{eq-jolie}
(b(x,y,\beta)\geq \alpha\;\mrm{and}\;\alpha,\beta>-\infty)
\iff 
(\bgal(y,x,\alpha)\geq \beta\;\mrm{and}\;\alpha,\beta>-\infty)\enspace.
\end{equation}
Moreover, 
\begin{equation}
\label{eq-jolie2}
(b(x,y,\beta)\geq \alpha\;\mrm{and}\;\alpha>-\infty)
\implies
\bgal(y,x,\alpha)\geq \beta\enspace. 
\end{equation}
\end{subequations}
\end{lemma}
\begin{proof}
We already observed in Remark~\ref{rem-sing}
that by~\eqref{a3} and~\eqref{e-bisinvertible}, $\bgal(y,x,\cdot)$
is the inverse of $b(x,y,\cdot)$
when $(x,y)\in\sS$.
Moreover, by Lemma~\ref{galoisr},
$b(x,y,\cdot)$ and $\bgal(y,x,\cdot)$ are nonincreasing.
When the left-hand side of~\eqref{eq-jolie2}
is satisfied, we have $(x,y)\in \sS$, 
so $\bgal(y,x,\cdot)$ is the inverse of $b(x,y,\cdot)$,
which shows~\eqref{eq-jolie2}. Together
with the symmetric implication, this shows~\eqref{eq-jolie}.
\end{proof}
\begin{proposition}\label{prop-partial}
If $g\in \sG$ and $g=BB\gal g$, then 
$(\partial\gal g )^{-1} = \partial B\gal g$.
\end{proposition}
\begin{proof} \sloppy
By~\eqref{e-inverse2}, 
$(\partial\gal g)^{-1}(y)=\set{x\in X}{(x,y)\in\sS\text{ and } 
\bgal(y,x,g(x))= B\gal g(y)}$. As noted in Lemma~\ref{galoisnondual},
for $(x,y)\in \sS$, $\bgal(y,x,\cdot)$
is a bijection with inverse $b(x,y,\cdot)$, so that
$(\partial\gal g)^{-1}(y)=\set{x\in X}{(x,y)\in\sS\text{ and } g(x)=
b(x,y,B\gal g(y))}$. Using $BB\gal g=g$ and~\eqref{e-inverse1},
we get $(\partial\gal g )^{-1}(y)=\partial B\gal g(y)$.
\end{proof}

\begin{remark}
When $B$ is the Legendre-Fenchel transform,
a function is in the image of $B$, or
of $B\gal$, 
if, and only if, it is either convex,
l.s.c., and proper, or identically $+\infty$, or identically
$-\infty$. Then Proposition~\ref{prop-partial}
gives the classical inversion
property of subdifferentials,
$(\partial g)^{-1}=\partial g^{\star}$,
which holds for all convex l.s.c.\ proper functions $g$
(see for instance~\cite[Theorem 23.5]{ROCK}).
\end{remark} 

\begin{proof}[Proof of \eqref{e-c0}$\Rightarrow$\eqref{e-eq1} of
Theorem~\ref{cover}]
Let us assume that $\partial\gal g(x)\neq\emptyset$,
for all $x\in X'$, and let us show
that $f=B\gal g$ satisfies $(\sP')$.
Since by~\eqref{galdef1},
$BB\gal g\leq g$, it is enough to  prove that
\begin{equation}
\partial\gal g (x)\neq\emptyset \impl BB\gal g(x) \geq  g(x)\enspace.
\label{tp3}
\end{equation}
If $y\in \partial\gal g(x)$, then, by~\eqref{e-inverse2},
$\bgal(y,x,g(x))=B\gal g(y)$ and $(x,y)\in\sS$, which yields
$BB\gal g(x) \geq b(x,y,B\gal g(y))=
b(x,y, \bgal(y,x,g(x)))=g(x)$,
by the first assertion of Lemma~\ref{galoisnondual},
and~\eqref{tp3} is shown.
\end{proof}

To pursue the proof of Theorem~\ref{cover}, 
we state properties of subdifferentials
which generalise Proposition~\ref{prop-partial}.
\begin{lemma}\label{lem-exis}
Consider $f\in \sF$ and $g\in \sG$ such that $Bf\leq g$, and let
$E=\set{x\in X}{Bf(x)=g(x)}$.
Then 
\begin{subequations}
\begin{align}
B B\gal g(x)&=g(x),\; \mrm{for all}\; x\in E,
\label{ttp2}
\end{align}
and for all $y\in Y$,
\begin{align}\partial f(y)\cap E&=
\begin{cases} \partial B\gal g(y)\cap E=(\partial\gal g)^{-1}(y)
& \mrm{if } f(y)= B\gal g(y),\\
\emptyset  &\mrm{otherwise.}
\end{cases}\label{ttp3}
\end{align}
\end{subequations}
\end{lemma}
\begin{proof}
By~\eqref{galdef2}, $Bf\leq g$ implies $f\geq B\gal g$.
Since $B$ is antitone,
applying $B$ to $f\geq B\gal g$, we get $Bf\leq B B\gal g$.
By~\eqref{galdef1}, $B B\gal g\leq g$,
hence $Bf\leq B B\gal g\leq g$, which implies~\eqref{ttp2}.

Let $x\in\partial f(y)\cap E$. Then, by~\eqref{e-inverse1},
$(x,y)\in \sS$ and $Bf(x)=b(x,y,f(y))$. Using the definition of $E$ and
the first assertion of Lemma~\ref{galoisnondual},
we get that  $f(y)=\bgal(y,x,g(x))$.
Using $f\geq B\gal g$ and the definition of $B\gal$, we obtain
$B\gal g(y)\leq f(y)=\bgal(y,x,g(x))\leq B\gal g(y)$.
Thus, $f(y)=B\gal g(y)$ and, using~\eqref{e-inverse2},
$y\in \partial\gal g(x)$.
This shows that $\partial f(y)\cap E\subset (\partial\gal g)^{-1}(y)$
for all $y\in Y$, and that $\partial f(y)\cap E=\emptyset$
when $f(y)\neq B\gal g(y)$.
To prove the converse inclusion in~\eqref{ttp3}, let $y\in Y$
such that $f(y)=B\gal g(y)$. Let $x\in (\partial\gal g)^{-1}(y)$,
then, by~\eqref{e-inverse2},
$(x,y)\in \sS$ and $f(y)=B\gal g(y)=\bgal(y,x,g(x))$.
It follows that $g(x)=b(x,y,f(y))\leq Bf(x)\leq g(x)$.
Hence, $x\in E$ and, by~\eqref{e-inverse1},  $x\in\partial f(y)$.
This shows that $\partial f(y)\cap E= (\partial\gal g)^{-1}(y)$
when $f(y)= B\gal g(y)$.
Replacing $f$ by $B\gal g$, we obtain that 
$\partial B\gal g(y)\cap \set{x\in X}{B B\gal g(x)=g(x)}=
 (\partial\gal g)^{-1}(y)$ for all $y\in Y$.
Taking the intersection with $E$, using~\eqref{ttp2} and using
$(\partial\gal g)^{-1}(y)=\partial f(y)\cap E\subset E$,
we obtain 
$\partial B\gal g(y)\cap E=(\partial\gal g)^{-1}(y)\cap E=
(\partial\gal g)^{-1}(y)$ when $f(y)= B\gal g(y)$.
\end{proof}

\begin{lemma}\label{parti}
Consider $f\in \sF$ and $g\in \sG$ such that $Bf\leq g$, and let
$E=\set{x\in X}{Bf(x)=g(x)}$. We have
\begin{subequations}
\begin{align}
& \left(\cup_{y\in f^{-1}(+\infty)} \partial f (y)\right)\cap E=
 \cup_{y\in (B\gal g)^{-1}(+\infty)} (\partial\gal g)^{-1}(y)
= g^{-1}(-\infty)\enspace,\label{e-incla}\\
&\left(\cup_{y\in \dom{f}} \partial f (y)\right)\cap E\subset
\cup_{y\in \dom{B\gal g}} (\partial\gal g)^{-1}(y)
\subset \dom{g}\enspace,\label{e-inclb}\\
&\left(\cup_{y\in f^{-1}(-\infty)} \partial f (y)\right)\cap E\subset
\cup_{y\in (B\gal g)^{-1}(-\infty)} (\partial\gal g)^{-1}(y)
\subset g^{-1}(+\infty)\enspace.\label{e-inclc}
\end{align}
\end{subequations}
\end{lemma}
\begin{proof}
The first inclusion in~\eqref{e-inclb} and \eqref{e-inclc},
and the inclusion of the left hand side term
of~\eqref{e-incla} in the middle term of~\eqref{e-incla},
follow readily from~\eqref{ttp3}.
If $x\in (\partial\gal g)^{-1}(y)$ with $y\in Y$, 
then, by~\eqref{e-inverse2},  $(x,y)\in\sS$ and 
$B\gal g(y) =\bgal(y,x,g(x))$.
Using that $\bgal(y,x,\cdot)$ is a decreasing
bijection $\rbar\to\rbar$, we get the second inclusions
in~\eqref{e-inclb} and \eqref{e-inclc}, 
and also the inclusion of the middle term of~\eqref{e-incla}
in the right hand side term of~\eqref{e-incla}.
This concludes the proof of~\eqref{e-inclb}
and \eqref{e-inclc}. It remains to prove the
inclusion of the right hand side term of~\eqref{e-incla}
in the left hand side term of~\eqref{e-incla}.
Consider $x\in X$ such that $g(x)=-\infty$. 
Since $Bf\leq g$, we get that $x\in E$.
Moreover, since $\sS_x\neq\emptyset$, 
there is a $y\in Y$ such that $b(x,y,\cdot)$
is bijective. This implies that
$f(y)\geq B\gal g(y)\geq \bgal(y,x,g(x))=
\bgal(y,x,-\infty)= +\infty$, hence $f(y)=+\infty$.
Moreover, $Bf(x)= g(x)=-\infty=b(x,y,+\infty)
=b(x,y,f(y))$, which shows, by~\eqref{e-inverse1}, that
$x\in \partial f(y)$. It follows that
$x\in \partial f(y)\cap E$ with $f(y)=+\infty$,
which concludes the proof.
\end{proof}

\begin{proof}[Proof of~\eqref{e-c0}$\Leftrightarrow$\eqref{e-c1}%
$\Leftrightarrow$\eqref{e-c2}$\Rightarrow$\eqref{e-c3}
and of~\eqref{e-c3}$\Rightarrow$\eqref{e-c2}, in Theorem~\ref{cover}.]\sloppy
The equivalence~\eqref{e-c0}$\Leftrightarrow$\eqref{e-c1} holds trivially
by the definition of $\dom{\partial\gal g}$ and of a covering.
By~\eqref{e-incla}, we get 
$\cup_{y\in (B\gal g)^{-1}(+\infty)} (\partial\gal g)^{-1}(y)
= g^{-1}(-\infty)=X\setminus \udom{g}$, and 
since $(B\gal g)^{-1}(+\infty)=Y\setminus \ldom{B\gal g}$,
we deduce \eqref{e-c1}$\Leftrightarrow$\eqref{e-c2}.
By~\eqref{e-inclc}, $\cup_{y\in (B\gal g)^{-1}(-\infty)} 
(\partial\gal g)^{-1}(y) \subset g^{-1}(+\infty)$, hence we always have
\eqref{e-c2}$\Rightarrow$\eqref{e-c3}. 
Since $\dom{B\gal g}\subset \ldom{B\gal g}$ and 
$X'\cap \udom g=X'\cap \dom g$ when $X'\subset \idom{g}\cup
g^{-1}(-\infty)$,
we have proved the 
implication~\eqref{e-c3}$\Rightarrow$\eqref{e-c2} 
in Theorem~\ref{cover}.
\end{proof}

Conditions~{\rm (\ref{e-eq1}--\ref{e-c3})} of Theorem~\ref{cover}
are trivial in the following degenerate cases:
\begin{proposition} Let $g\in\sG$. We have 
\begin{subequations}
\begin{align}
\label{a2} 
g\equiv +\infty \;\Leftrightarrow\; B\gal g\equiv -\infty\enspace,\\
\label{a4} g\equiv -\infty \impl B\gal g\equiv +\infty\enspace.
\end{align}
\end{subequations}
In both cases,
 $BB\gal g=g$ and
$\{(\partial\gal g)^{-1}(y)\}_{y\in Y}$ is a covering of $X$.
Moreover,
if $B\gal g\equiv +\infty $ and $BB\gal g=g$,
then $g\equiv -\infty$.
\end{proposition}
\begin{proof} 
The implication $\Rightarrow$ in~\eqref{a2} follows from~\eqref{e-toptobot}.
By symmetry, if $f\in\sF$ then $f\equiv +\infty$ implies $Bf\equiv-\infty$.
Taking $f=B\gal g$, we get 
\begin{equation}
B\gal g\equiv+\infty \implies BB\gal g\equiv -\infty\enspace,
\label{e-util}
\end{equation}
which implies the last assertion of the lemma.
If $g\equiv -\infty$, then for all $y\in Y$,
taking $x\in \sS^y$, we get 
$B\gal g(y)\geq \bgal(y,x,g(x))=+\infty$,
which shows~\eqref{a4}.
By symmetry, if $f\in\sF$ then $f\equiv -\infty$ implies
$Bf\equiv +\infty$. Applying
this property to $f=B\gal g$, we get
that $B\gal g\equiv -\infty$ implies
$BB\gal g\equiv +\infty$, and since
$g\geq BB\gal g$, $g\equiv +\infty$,
which shows the implication $\Leftarrow$
in~\eqref{a2}, together with $BB\gal g=g$.
When $g\equiv -\infty$, combining~\eqref{a4}
and~\eqref{e-util}, we also get $BB\gal g=g$.
Moreover, since $\udom g=\emptyset$, \eqref{e-c2} is trivial with $X'=X$, and
by the equivalence \eqref{e-c1}$\Leftrightarrow$\eqref{e-c2}, which
has already been proved, we get that 
$\{(\partial\gal g)^{-1}(y)\}_{y\in Y}$ is a covering of $X$.
It remains to show that the same
covering property holds when $g\equiv +\infty$.
For all $x\in X$ and $y\in Y$, 
we have $B\gal g(y)= -\infty=
\bgal(y,x,g(x))$. Taking $y\in \sS_x$,
we get $y\in \partial\gal g(x)$ by~\eqref{e-inverse2},
which shows that
$\cup_{y\in Y} (\partial\gal g)^{-1}(y)=X$
in this case, too.
\end{proof}

We next mention some direct consequences
of the continuity and coercivity assumptions.
\begin{lemma}\label{lem-usctocon}
The kernel $b$ is continuous in the second variable if, and
only if, for all $x\in X$ and $\alpha\in \R$,
$b(x,\cdot,\alpha)$ is upper semicontinuous (u.s.c.).
In that case, $b(x,\cdot,\alpha)$ is continuous,
for all $x\in X$ and $\alpha\in \R\cup\{+\infty\}$.
\end{lemma}
\begin{proof}
Proposition~\ref{prop-blsc}
shows that $b(x,\cdot,\alpha)$ 
is l.s.c.\ for all $x\in X$ and $\alpha\in \rbar$.
Hence, for all $x\in X$ and $\alpha\in \R$,
$b(x,\cdot,\alpha)$ is u.s.c.\ if, and only if, it is
continuous.
Moreover, since, by~\eqref{e-toptobot}, 
$b(x,\cdot,+\infty)\equiv -\infty$, 
$b(x,\cdot,+\infty)$ is always (trivially) continuous.
\end{proof}
Note that the continuity assumption
does {\em not} require that
$b(x,\cdot,-\infty)$ is continuous or u.s.c.\ 
(indeed, in the special case
when $b(x,\cdot,\alpha)=
\bar b(x,\cdot)-\alpha$,
we have $b(x,\cdot,-\infty)=
\bar b(x,\cdot)+\infty$, which need
not be u.s.c.\ if $\bar b(x,\cdot)$ is continuous
and takes the value $-\infty$).
The next lemma shows that 
assuming $b$ or $\bgal$ to be u.s.c.\ (or equivalently continuous) is the same:
\begin{lemma}\label{lem-dualusc}
Let $x\in X$. Then $b(x,\cdot,\beta)$ is u.s.c.\ 
for all $\beta\in \R$
if, and only if, $\bgal(\cdot,x,\alpha)$ is u.s.c.\ 
for all $\alpha\in \R\,$.
\end{lemma}
\begin{proof}
\sloppy
We already observed in the proof of Lemma~\ref{lem-usctocon}
that $b(x,\cdot,+\infty)$ and $\bgal(\cdot,x,+\infty)$ are
u.s.c., so it is enough to show that
\begin{align}
\label{e-usctousc}
(b(x,\cdot,\beta) \;\mrm{is u.s.c.}\;
\forall \beta\in \R\cup\{+\infty\} )
\iff 
(\bgal(\cdot,x,\alpha) \;\mrm{is u.s.c.}\;
\forall \alpha\in \R\cup\{+\infty\} )\enspace.
\end{align}
The left hand side of~\eqref{e-usctousc}
is equivalent to 
\begin{equation}
\label{e-usc}
\set{y\in Y}{b(x,y,\beta)\geq \alpha}\;\mrm{is closed}\;
\forall \alpha,\beta\in \R\cup\{+\infty\}\enspace.
\end{equation}
Applying~\eqref{eq-jolie},
we get that~\eqref{e-usc} is equivalent to
\[
\set{y\in Y}{\bgal(y,x,\alpha)\geq \beta}\;\mrm{is closed}\;
\forall \alpha,\beta \in \R\cup\{+\infty\}\enspace,
\]
which is exactly
the upper semicontinuity of all the maps
$\bgal(\cdot,x,\alpha)$, with $\alpha \in \R\cup\{+\infty\}$.
\end{proof}
\begin{lemma}\label{lem-compose}
If $b$ is continuous in the second variable, then for all maps $f\in \sF$
such that $f(y)>-\infty$ for all $y\in Y$, 
the map $y\mapsto b(x,y,f(y))$ is u.s.c.\ for all $x\in X$.
\end{lemma}
\begin{proof}
We have to show that $\set{y\in Y}{b(x,y,f(y))\geq \beta}$
is closed, for all $\beta \in \R\cup\{+\infty\}$.
Since $f(y)>-\infty$, for all $y\in Y$,~\eqref{eq-jolie} yields
$\set{y\in Y}{b(x,y,f(y))\geq \beta}=
\set{y\in Y}{\bgal(y,x,\beta) \geq f(y)}$,
which is closed since $f$ is l.s.c.\ and 
$\bgal(\cdot,x,\beta)$ is u.s.c.\ (by Lemma~\ref{lem-usctocon}).
\end{proof}
The following observation shows that we
could have replaced ``relatively compact''
by ``compact'' in the definition of
coercivity.
\begin{proposition}
If $b$ is continuous in the second variable and coercive, 
then, for all $\alpha,\beta\in \R$,
for all $x\in X$ and neighbourhoods $V$ of $x$,
$\set{y\in Y}{b_{x,V}^{\alpha}(y)\leq \beta}$
is compact.
\end{proposition}
\begin{proof}
Since $b_{x,V}^{\alpha}$ is given by the sup in~\eqref{e-def-bv},
we have $\set{y\in Y}{b_{x,V}^{\alpha}(y)\leq \beta}=
\cap_{z\in V} Y_z$,
where $Y_z=\set{y\in Y}{ b(z,y,\bgal(y,x,\alpha))\leq \beta}$.
By~\eqref{galdef2},
 $Y_z=\set{y\in Y}{\bgal(y,z,\beta) \leq \bgal(y,x,\alpha)}$,
which is closed because $\bgal(\cdot,z,\beta)$ is 
l.s.c.\ (by Theorem~\ref{theo-rep}), and
$\bgal(\cdot,x,\alpha)$ is u.s.c.\ for $\alpha\in \R$
(by Lemma~\ref{lem-dualusc} and the continuity of $b$ in the second variable).
Therefore, $\cap_{z\in V} Y_z$,
which is closed and relatively compact, is compact.
\end{proof}

The proof of \eqref{e-eq1}$\Rightarrow$\eqref{e-c2} in Theorem~\ref{cover}
relies on the following result:
\begin{theorem}\label{exis}
Let $f\in\sF$.
Assume that either $Y$ is discrete, or $b$ is continuous in the second variable
and $f(y)>-\infty$ for all $y\in Y$.
Then, if $f\in\sFC$, 
$\{\partial f(y)\}_{y\in \ldom{f}}$ is a covering of $\udom {Bf}$,
and if $b$ is coercive, 
$\{\partial f(y)\}_{y\in \ldom{f}}$ is a covering of $\idom {Bf}$.
\end{theorem}
\begin{proof}
We set $g=Bf$. We prove at the same time the two assertions
of the theorem by setting $X_0=\udom g$ when $f\in\sFC$, and $X_0=\idom g$ 
when $b$ is coercive. We thus need to prove that 
$X_0\subset \cup_{y\in\ldom{f} }\partial f(y)$.
Since, by~\eqref{e-incla},  $\cup_{y\in Y\setminus\ldom{f} }\partial f(y)
=X\setminus \udom g$, and since $X_0\subset  \udom g$, it is sufficient
to prove that $X_0\subset \cup_{y\in Y}\partial f(y)$.
We will prove:
\begin{equation}
x\in X_0 \implies \exists y\in Y,\;  b(x,y,f(y))= 
\sup_{y'\in Y} b(x,y',f(y'))\enspace.
\label{e-toprove}
\end{equation}
Indeed, if~\eqref{e-toprove} holds, then for all $x\in X_0$, there 
exists $y\in Y$ such that $b(x,y,f(y))=Bf(x)=g(x)$ and
since $X_0\subset \udom g$, $g(x)\neq -\infty$.
Hence $b(x,y,f(y))\neq -\infty$, whence
$(x,y)\in \sS$, which implies with $b(x,y,f(y))=Bf(x)$
that $x\in \partial f(y)$. This shows that 
$X_0\subset \cup_{y\in Y}\partial f(y)$.

To prove~\eqref{e-toprove}, it suffices
to show that 
\begin{equation}
\forall x\in X_0,\;\forall \alpha\in \R,\;\;
L_{\alpha}(x)=\set{y\in Y}{b(x,y,f(y))\geq \alpha}
\;\mrm{is compact.}
\label{e-toprove2}
\end{equation}
Let us first prove that the sets $L_{\alpha}(x)$ are closed 
for all $x\in X$ and $\alpha\in\R$.
When $Y$ is discrete, this is trivial.
Otherwise, by the assumptions of the theorem, 
$f(y)>-\infty$, for all $y\in Y$, 
and $b$ is continuous in the second variable, therefore, by
Lemma~\ref{lem-compose}, $y\mapsto b(x,y,f(y))$ is an u.s.c.\ map
for all $x\in X$. This implies again that the sets $L_{\alpha}(x)$ are closed 
for all $x\in X$ and $\alpha\in\R$.

It remains to show that the sets $L_{\alpha}(x)$ are relatively
compact for all $x\in X_0$ and $\alpha\in\R$. 
By definition of $\sFC$, this holds trivially for
any $X_0\subset X$, when $f\in\sFC$. Let us finally
assume that $b$ is coercive and $X_0=\idom g$.
Let $x\in\idom g$ and $\alpha\in\R$. There exists $\beta\in \R$ 
such that $\limsup_{x'\to x}  g(x') <\beta$, so there exists 
a neighbourhood $V$ of $x$ in $X$ such that
 $\sup_{x'\in V}g(x')\leq \beta$.
Then, by~\eqref{eq-jolie2}:
\[b(x,y,f(y)) \geq \alpha \implies \bgal(y,x,\alpha) \geq f(y)\enspace,\]
and since
\[ 
f(y)\geq B\gal g (y) \geq \bgal(y,z,g(z)) 
\geq \bgal(y,z,\beta)\quad \forall z\in V\enspace,
\]
we obtain:
\begin{align*} b(x,y,f(y)) \geq \alpha &  \implies 
\forall z\in V,\; \bgal(y,x,\alpha)\geq \bgal(y,z,\beta) \\
&\implies \forall z\in V,\;
\beta\geq b(z,y,\bgal(y,x,\alpha))\enspace,
\end{align*}
which shows that 
$L_\alpha(x)\subset \set{y\in Y}{b_{x,V}^{\alpha}(y)\leq \beta}$.
By the coercivity of $b$, the latter set is relatively compact,
and thus $L_\alpha(x)$ is also relatively compact.
This concludes the proof of~\eqref{e-toprove2}.
\end{proof}

\begin{proof}[Proof of \eqref{e-eq1}$\Rightarrow$\eqref{e-c2} 
in Theorem~\ref{cover}]\sloppy
If $(\sP')$  has a solution $f\in \sF$, then, by Lemma~\ref{lem-exis}
and~\eqref{galdef1}, 
$f=B\gal g$ is also a solution of $(\sP')$. Fix $f=B\gal g$.
By Lemma~\ref{parti},
$\left(\cup_{y\in \ldom{f}}\partial f(y)\right)\cap X'
\subset \cup_{y\in \ldom{B\gal g}}(\partial\gal g)^{-1}(y)$.

Let us first consider the case where Assumption~\ref{casei} holds.
Then, applying Theorem~\ref{exis} to $f$,
we deduce that 
 $\udom{Bf}\cap X'\subset\cup_{y\in \ldom{B\gal g}}(\partial\gal g)^{-1}(y)$.
Since it is clear that $X'\cap \udom g\subset X'\cap \udom {Bf}$,
we get~\eqref{e-c2}.

Let us finally consider the case where Assumption~\ref{caseii} holds.
Then, applying Theorem~\ref{exis} to $f$,
we obtain that
$\idom{Bf}\cap X'\subset \left(\cup_{y\in \ldom{f}} \partial f(y)\right)
\cap X'\subset\cup_{y\in \ldom{B\gal g}}(\partial\gal g)^{-1}(y)$.
Moreover, it is easy to show that $X'\subset \idom{g}\cup g^{-1}(-\infty)$
implies $X'\cap \udom {g} \subset X'\cap \idom{Bf}$,
hence~\eqref{e-c2} follows.
\end{proof}

\begin{proof}[Proof of Corollary~\ref{cor-compact}]
When Assumptions~
~\ref{a-bgalgisfinite}, and \ref{a-bgalginfc}, hold, 
the equivalence \eqref{e-eq1}$\Leftrightarrow$\eqref{e-c1} 
in Theorem~\ref{cover}, for $X'=X$,
yields the first assertion of the corollary.
Assume now that $\{( \partial\gal g)^{-1}(y)\}_{y\in Y}$
is a covering of $X$. Then, for all $x\in X$,  there exists
$y\in \partial\gal g (x)$. By~\eqref{e-inverse2},
 $(x,y)\in\sS$ and $B\gal g(y)= \bgal(y,x,g(x))$.
Since $B\gal g(y)>-\infty$ for all $y\in Y$, and 
$\bgal(y,x,\cdot)$ is a decreasing bijection $\rbar\to\rbar$,
we deduce that $g(x)<+\infty$, for all $x\in X$.
\end{proof}

\section{Uniqueness of Solutions of $Bf=g$}\label{uni-sec}
\label{sec-unique}
\subsection{Statement of the Uniqueness Results}
To give a uniqueness result,
we need some additional definitions.

\begin{definition}\label{defi-cover2}
Let $F$ be a map from a set $Z$ to the set $\pW$ of all subsets of some
set $W$, and let $Z'\subset Z$ and $W'\subset W$ be such that 
$\{F(z)\}_{z \in Z'}$ is a covering of $W'$.
An element $y\in Z'$ is said {\em algebraically essential} with respect to
the covering $\{F(z)\}_{z \in Z'}$ of $W'$
if there exists $w\in W'$ such that
$w\not\in \cup_{z\in Z'\setminus\{y\}} F(z)$.
When $Z$ is a topological space,
an element $y\in Z'$ is said {\em topologically essential} with respect to
the covering $\{F(z)\}_{z \in Z'}$ of $W'$ if for all open neighbourhoods
$U$ of $y$ in $Z'$, there exists $w\in W'$ such that
$w\not\in \cup_{z\in Z'\setminus U} F(z)$.
The covering of $W'$ by $\{F(z)\}_{z \in Z'}$ is 
{\em algebraically (resp.\ topologically) minimal} if all elements
of $Z'$ are algebraically (resp.\ topologically) essential.
\end{definition}
Algebraic minimality implies topological minimality.
Both notions coincide if $Z$ is a discrete topological space.

\begin{definition}\label{expos}
Let $f\in\sF$ and $X'\subset X$.
We say that $y\in Y$ is an \NEW{exposed point} of $f$ relative to $X'$
if there exists $x\in X'$ such that $(x,y)\in \sS$ and
\[
 b(x,y',f(y'))< b(x,y,f(y)) \quad \forall y'\in Y\setminus\{y\}\enspace.
\]
\end{definition}
When $B$ is the Legendre-Fenchel transform and $X'=X$, this notion coincides 
with the definition given in~\cite[Definition 2.3.3]{DEMBO}
of an exposed point of $f$. It is equivalent to the property that
$(y,f(y))$ is an exposed point of the epigraph of 
$f$~\cite[Sections 18 and 25]{ROCK}.
We readily get from Definitions~\ref{defi-cover} and~\ref{expos}:
\begin{lemma}
Let $f\in\sF$ and let $X'\subset \cup_{y\in Y} \partial f(y)$.
An element $z\in Y$ is an exposed point of $f$ relative to $X'$
if, and only if, $z$ is algebraically essential with respect to
the covering $\{\partial f(y)\}_{y \in Y}$ of $X'$.
\end{lemma}

\begin{definition}
Let $Z$ and $W$ be topological spaces.
We say that a map $h:Z\to W$ is \NEW{quasi-continuous}
if for all open sets $G$ of $W$, the set $h^{-1}(G)$ 
is \NEW{semi-open},
which means that $h^{-1}(G)$ is included in the closure of its interior.
\end{definition}
See for instance \cite{NEUBRUNN} for definitions and 
properties of quasi-continuous functions or multi-applications.
If $h:X\to\rbar$ is l.s.c., then $h$ is quasi-continuous if, and only if,
$h=\lsc(\usc(h))$, where $\lsc$ (resp.\ $\usc$)
 means the l.s.c.\ (resp.\ u.s.c.) hull.
The notion of quasi-continuous function, and the properties of l.s.c.\ or
u.s.c.\ quasi-continuous functions have also  been
studied in~\cite{samb02}.
\begin{definition}
We say that $B$ is \NEW{regular} if
for all $f\in\sF$, $B f$ is  l.s.c.\ on $X$ and
 quasi-continuous on its domain,
which means that the restriction
of $Bf$ to its domain is quasi-continuous for the induced topology.
\end{definition}
The notion of regularity
for $B\gal$ is defined in the symmetric way.
When $X$ (resp.\ $Y$) is endowed with the discrete topology,
$B$ (resp.\ $B\gal$) is always regular.
When $\sS=X\times Y$ and $\{b(\cdot,y,\alpha)\}_{y\in Y,\; \alpha\in\R}$ 
is an equicontinuous family of functions, then 
$Bf$ is continuous on $X$ for any $f\in\sF$, so $B$ is regular.
The Legendre-Fenchel transform on $\R^n$ is regular (see Lemma~\ref{fenchel-reg}
below).

We now state several uniqueness results, that we shall prove
in Section~\ref{sec-proof-unique}. 
\begin{theorem}\label{minicover1}
Let $X'\subset X$, and $g\in \sG$.
Assume that $\{ (\partial\gal g)^{-1}(y)\}_{y\in \ldom{B\gal g}}$ 
is a covering of $X' \cap \udom g$, and denote by
$Z_a$ (resp.\ $Z_t$) the set of algebraically (resp.\ topologically)
essential elements with respect to this covering.
Make Assumptions~
\ref{a-yisdiscrete} or \ref{a-bgalgisfinite},
and \ref{a-bgalginfc} or \ref{a-biscoercive}.
Then Problem  $(\sP')$ has a solution, and 
any solution  $f$ of $(\sP')$ satisfies
\begin{align} \label{psec}
f\geq B\gal g,\quad\mrm{and}\quad f(y)=B\gal g(y) \quad
\mrm{ for all }y\in Z_a\enspace.
\end{align}
If, in addition, $B\gal g$ is quasi-continuous on its domain,
and $\intt{Z_t}$ denotes the interior of $Z_t$, relatively to
$\ldom{B\gal g}$, then any solution  $f$ of $(\sP')$ satisfies
\begin{align} \label{psec2}
f(y)=B\gal g(y) \quad\text{ for all }y\in \intt{Z_t}\enspace.
\end{align}
\end{theorem}
\begin{theorem}\label{minicover}
Let $X'\subset X$, and $g\in \sG$.
Consider the following statements:
\begin{align}
&\mrm{Problem }(\sP')\mrm{ has a unique solution}, \label{fe-eq1} \\
&\{ (\partial\gal g)^{-1}(y)\}_{y\in \ldom{B\gal g}} 
\mrm{ is a topologically minimal covering of }\label{fe-c1} \\ &
 X' \cap \udom g. \nonumber
\end{align}
We have {\rm (\ref{fe-eq1},\ref{e-c2})}$\Rightarrow$\eqref{fe-c1}.
The implication \eqref{fe-eq1}$\Rightarrow$\eqref{fe-c1}
holds if Assumptions~
\ref{a-yisdiscrete} or \ref{a-bgalgisfinite},
and \ref{a-bgalginfc} or \ref{a-biscoercive},
are satisfied.
The equivalence \eqref{fe-eq1}$\Leftrightarrow$\eqref{fe-c1}
holds if we assume in addition that 
$B\gal g$ is quasi-continuous on its domain.  In particular,
this equivalence holds when $Y$ is finite.
\end{theorem}

The topological minimality in~\eqref{fe-c1} 
is a relaxation of algebraic minimality,
which is a generalised differentiability condition.
Indeed, if $\{ (\partial\gal g)^{-1}(y)\}_{y\in \ldom{B\gal g}}$
is a covering of $X' \cap \udom g$, this covering
is algebraically minimal
if, and only if, for all $y\in \ldom{B\gal g}$, 
there is an $x\in X'\cap \udom g$ such that
$\partial\gal g(x)=\{y\}$.
This is in particular fulfilled when
$\{ (\partial\gal g)(x)\}_{x\in X'\cap\udom{g}}$ is a covering of 
$\ldom{B\gal g}$, and for all $x\in X'\cap \udom{g}$,
$\partial\gal g(x)$ is a singleton, a 
condition which, in the case where $B$ is the Legendre-Fenchel transform,
means that $g$ is differentiable at $x$.

\begin{corollary}\label{unifini}
Consider $g\in \sG$. Assume that $Y$ is finite.
Then the equation $Bf=g$ has a unique solution 
$f\in\sF$, if, and only if, 
$\{(\partial\gal g)^{-1}(y)\}_{y\in \ldom{B\gal g}}$ is an
algebraically minimal covering of $\udom g$.
\end{corollary}

\begin{corollary}\label{uniccompact}
Let  $g\in \sG$. Make Assumptions
~\ref{a-bgalgisfinite}, and \ref{a-bgalginfc}.
Assume in addition that $B\gal g$ is quasi-continuous on its domain.
Then the equation $Bf=g$ has a unique solution 
$f\in\sF$, if, and only if, 
$\{(\partial\gal g)^{-1}(y)\}_{y\in \ldom{B\gal g}}$ is a 
topologically minimal covering of $\udom g$.
\end{corollary}

Since \eqref{fe-eq1} implies that
Problem $(\sP)$ has at most one solution,
Theorem~\ref{minicover} yields a sufficient
condition for the uniqueness of the solution of Problem $(\sP)$.
However, for Problem $(\sP)$, the necessary uniqueness condition
implied by Theorem~\ref{minicover} only holds when $B\gal g\in\sFC$, or when 
$X=\idom{g}\cup g^{-1}(-\infty)$, or when $X=\dom{\partial\gal g}$.
To give a more specific uniqueness result for Problem $(\sP)$,
we shall use the following condition:
\[
\condcmath:\quad 
\begin{array}{l}
\mrm{there exists a basis }\sB\mrm{ of neighborhoods such that}\\[1mm]
\displaystyle 
\forall U\in\sB,\; \exists \varepsilon>0,\;
\forall x\in X,\;
\sup_{y\in U\cap \sS_x,\, \alpha\in\R} (b(x,y,\alpha)-
b(x,y,\alpha+\varepsilon))<+\infty\enspace.
\end{array}
\]
Condition \condc\ is satisfied in particular when
$\{b(x,y,\cdot)\}_{x\in X,\; y\in Y}$ is a family of $\beta$-H\"older 
continuous functions (for $0<\beta\leq 1$), uniformly in $y\in U$, 
for all small 
enough open sets $U$, or if $\{b(x,y,\cdot)\}_{x\in X,\; y\in U}$ 
is an equicontinuous family, for all small enough open sets $U$.
In particular, condition \condc\ is satisfied when 
$b(x,y,\alpha)=b(x,y)-\alpha$ or when 
$b(x,y,\alpha)=-(|\<x,y>|+1)\alpha$.

\begin{theorem}\label{theo-uniquenew}
Let $g\in\sG$. Then the existence and uniqueness of a solution
of Problem $(\sP)$ implies \eqref{fe-c1},
if one of the following assertions is satisfied:
\begin{enumerate}
\renewcommand{\theenumi}{\arabic{enumi}}
\item\label{lodiscrete}
$X'=\dom{\partial\gal g}$, and Assumption~\ref{a-yisdiscrete} holds;
\item\label{locompact}
$X'=\dom{\partial\gal g}$, $Y$ is locally compact,
and Assumption
~\ref{a-bgalgisfinite} holds;
\item\label{regular} $X'=\idom g$,
$b$ is coercive, $B$ is regular,
Assumption
~\ref{a-bgalgisfinite} holds,
$\dom g$ is included in the closure of $\idom g$,
and either $Y$ is locally compact or
condition \condc\ holds.
\end{enumerate}
\end{theorem}

\subsection{Proofs of the Uniqueness Results}\label{sec-proof-unique}
Let us first state a general property of quasi-continuous functions. 
\begin{lemma}\label{semiopen}
Let $f\in \lsc(Y,\rbar)$, and let
$h:Y\to\rbar$ be a quasi-continuous map. Then the set
$V=\set{y\in Y}{h(y)< f(y)}$ is semi-open.
In particular, if $V$ is non-empty, $V$ has a non-empty
interior.
\end{lemma}
\begin{proof}
It suffices to consider the case when $V$ is non-empty.
Let $z\in V$.
There exists $a\in\R$ such that $h(z)<a<f(z)$.
Consider $V_1=\set{y\in Y}{h(y)<a}$,
$U_1$ the interior of $V_1$,
and $U_2=\set{y\in Y}{a<f(y)}$.
We have $z\in V_1\cap U_2\subset V$.
Since $h$ is quasi-continuous, $V_1$ is semi-open, hence
$V_1$ is included in the closure
of $U_1$, that we denote by $\overline{U_1}$.
Since $f$ is l.s.c., $U_2$ is open.
We have $z\in V_1\cap U_2
\subset \overline{U_1}\cap U_2
\subset \overline{U_1\cap U_2}$,
and since $U_1\cap U_2$ is open
and included in $V$, we get that
$z$ belongs to the closure of the interior of $V$.
\end{proof}

\begin{proof}[Proof of Theorem~\ref{minicover1}]
Note first that since, when $Bf \leq g$, the equation $Bf=g$ on $X'$ is
equivalent to the  equation $Bf=g$ on $X'\cap \udom g$,
we can assume without restriction of generality that
$X'\subset \udom g$ in the proofs of Theorems~\ref{minicover1}
and~\ref{minicover}.

Let $Y'=\ldom{B\gal g}$, 
assume that $\{ (\partial\gal g)^{-1}(y)\}_{y\in Y'}$ 
is a covering of $X'$, and denote by
$Z_a$ (resp.\ $Z_t$) the set of algebraically (resp.\ topologically)
essential elements with respect to this covering.
Since~\eqref{e-c2} holds, it follows from Theorem~\ref{cover}, 
that $(\sP')$ has a solution $f\in \sF$, 
and by Lemma~\ref{lem-exis} and~\eqref{galdef1}, $B\gal g$
is necessarily  another solution and $B\gal g\leq f$.

Let $f\in \sF$ be a solution of $(\sP')$ and denote by
$F=\set{y\in Y'}{B\gal g(y)=f(y)}$ and by $V$ its complement in $Y'$.
Since, $B\gal g\leq f$, $V=\set{y\in Y'}{B\gal g(y)<f(y)}$.
We claim that 
\begin{equation}\label{minieq1}
\forall x\in X' \; \exists y\in F\mrm{ such that }
x\in (\partial\gal g)^{-1} (y)\enspace.
\end{equation}
If~\eqref{minieq1} is proved, then the following holds
\begin{equation}\label{minieq2}
Z_a\subset F\mrm{ and } Z_t\subset \overline{F}\enspace,
\end{equation}
where $\overline{F}$ denotes the closure of $F$, relatively to $Y'$.
Indeed, let us first consider $z\in Z_a$. Since $z$ is algebraically essential
with respect to the covering 
$\{(\partial\gal g)^{-1}(y)\}_{y\in Y'}$  of $X'$,
there exists $x\in X'$  such that
\begin{align}\label{unic0} x\in 
(\partial\gal g)^{-1}(z)\setminus \cup_{y\in Y'\setminus \{z\}}
(\partial\gal g)^{-1}(y)\enspace. \end{align}
Moreover, by~\eqref{minieq1}, there exists $y\in F$ such that
$x\in (\partial\gal g)^{-1} (y)$. Using~\eqref{unic0}, this implies
that $y=z$ and $z\in F$, which shows $Z_a\subset F$.
Let us now consider $z\in Z_t$. Since $z$ is topologically essential
with respect to the covering 
$\{(\partial\gal g)^{-1}(y)\}_{y\in Y'}$  of $X'$, for
all open neighbourhoods $U$ of $z$ in $Y'$, 
there exists $x_U\in X'$  such that
\begin{align}\label{unic1} x_U\in \cup_{y\in U} 
(\partial\gal g)^{-1}(y)\setminus \cup_{y'\in Y'\setminus U}
(\partial\gal g)^{-1}(y')\enspace. \end{align}
Moreover, by~\eqref{minieq1}, there exists $y_U\in F$ such that
$x_U\in (\partial\gal g)^{-1} (y_U)$. Using~\eqref{unic1}, this implies
that $y_U\in U$. We have thus
proved that for all open neighbourhoods $U$ of $z$ in $Y'$,
there exists $y_U\in U\cap F$, which means that $z\in \overline{F}$,
and shows $Z_t\subset \overline{F}$.

Now from~\eqref{minieq2}, we get~\eqref{psec}.
If Assumption~\ref{a-yisdiscrete} holds,
that is, if $Y$ is discrete, $Z_a=Z_t$ and thus~\eqref{psec2} holds trivially
by~\eqref{psec}.
Otherwise, we deduce from Assumption~\ref{a-bgalgisfinite}
that $Y'=\ldom{B\gal g}=\dom{B\gal g}$. Moreover, if
$B\gal g$ is quasi-continuous on its domain,
then, since $f$ is l.s.c.\ on $Y$, and a 
fortiori on $Y'$, we obtain, by Lemma~\ref{semiopen},
that $V$ is semi-open in $Y'$.
It follows that its complement in $Y'$, $F$, contains the interior of its 
closure $\overline{F}$, relatively to $Y'$.
Using~\eqref{minieq2}, this yields $\intt{Z_t}\subset F$,
which means precisely that~\eqref{psec2} holds.

Let us prove~\eqref{minieq1}.
When $f$ is a solution of $(\sP')$, $X'\subset \udom g$ implies
$X'\subset \udom {Bf}$, $X'\subset \idom g$ implies
$X'\subset \idom{Bf}$, and $B\gal g\in\sFC$ implies
$f\in\sFC$ (since $f\geq B\gal g$). Hence, Theorem~\ref{exis} shows that
$X'\subset \cup_{y\in \ldom{f}} \partial f(y)$.
Since $B\gal g\leq f$ implies that $\ldom{f}\subset \ldom{B\gal g}=Y'$,
we get that $X'\subset \cup_{y\in Y'} \partial f(y)$.
Hence, for all $x\in X'$, there exists $y\in Y'$ such that
$x\in\partial f(y)$.
Since then $\partial f(y)\cap X'\neq\emptyset$, Lemma~\ref{lem-exis}
shows that $f(y)=B\gal g(y)$
 and $\partial f(y)\cap X'\subset (\partial\gal g)^{-1}(y)$.
Hence, $y\in F$ and $x\in (\partial\gal g)^{-1}(y)$,
which shows~\eqref{minieq1}.
\end{proof}

We now prove the different assertions of Theorem~\ref{minicover}.

\begin{proof}[Proof of {\rm (\ref{fe-eq1},\ref{e-c2})}$\Rightarrow$\eqref{fe-c1} in Theorem~\ref{minicover}]
We assume, as in the above proof, that $X'\subset  \udom g$.
Set $Y'=\ldom {B\gal g}$.
Assume that~\eqref{fe-eq1} and~\eqref{e-c2} hold, 
which means that $(\sP')$ has a unique solution and 
$\{(\partial\gal g)^{-1}(y)\}_{y\in Y'}$ is a covering of $X'$.
Assume by contradiction that this covering
is not topologically minimal, 
i.e., that there exists an open set $U$ of $Y$ such that
$U\cap Y'\neq\emptyset$, and 
such that for all $x\in X'$,
there exists $y\in Y'\setminus U$ such that
$x\in (\partial\gal g)^{-1}(y)$,  which means, by~\eqref{e-inverse2}, 
that $(x,y)\in\sS$ and $B\gal g(y)=\bgal(y,x,g(x))$.
Then  $g(x)=b(x,y,B\gal g(y))$ and, since $g\geq B B\gal g$,
we get:
\begin{equation}
g(x)=\sup_{y\in Y\setminus U} b(x,y,B\gal g(y))
\quad\forall x\in X'\enspace.
\label{eq-n}
\end{equation}
To contradict the uniqueness
of the solution of $(\sP')$, it suffices to construct a map 
$f\in\sF$ such that $f\neq B\gal g$ and 
\begin{equation}
f\geq B\gal g, \qquad f=B\gal g \;\mrm{on}\; Y\setminus U
\label{keyprop}\enspace.
\end{equation}
Indeed, for any function $f$ satisfying~\eqref{keyprop},
we have $Bf\leq g$, and, by~\eqref{eq-n}, 
$Bf\geq g$ on $X'$, hence $f$ is a solution of $(\sP')$.
The function $f=B\gal g$ satisfies trivially $f\in\sF$ and \eqref{keyprop}.
Defining $f$ by $f=B\gal g$ on $Y\setminus U$ and $f=+\infty$ on $U$,
we obtain that $f\in\sF$, $f$ satisfies \eqref{keyprop} and since
$Y'\cap U\neq \emptyset$, $f\neq B\gal g$, which concludes the proof.
\end{proof}

\begin{proof}[Proof of \eqref{fe-eq1}$\Rightarrow$\eqref{fe-c1}
in Theorem~\ref{minicover}]
If \eqref{fe-eq1} holds, we get in particular that Problem~$(\sP')$
has a  solution, 
and we deduce~\eqref{e-c2} from Theorem~\ref{cover},
using the Assumptions of Theorem~\ref{minicover}.
{From} the implication {\rm (\ref{fe-eq1},\ref{e-c2})}$\Rightarrow$%
\eqref{fe-c1}  that we already proved, we obtain 
\eqref{fe-eq1}$\Rightarrow$\eqref{fe-c1}.
\end{proof}

\begin{proof}[Proof of \eqref{fe-c1}$\Rightarrow$\eqref{fe-eq1}
in Theorem~\ref{minicover}]
The assumptions for this implication, in Theorem~\ref{minicover}, imply
that the assumptions of Theorem~\ref{minicover1} are satisfied.
They also imply that
any  element of $\ldom{B\gal  g}$ is topological essential 
for the covering $\{(\partial\gal g)^{-1}(y)\}_{y\in \ldom{B\gal  g}}$
of $X'\cap\udom{g}$. This means that $Z_t=\ldom{B\gal  g}$ in
Theorem~\ref{minicover1}, thus
$\intt{Z_t}=Z_t=\ldom{B\gal g}$,
hence  the conclusions of Theorem~\ref{minicover1}
imply that $(\sP')$ has a solution
and that any solution $f$ of $(\sP')$ satisfies
$f\geq B\gal g$ and $f(y)=B\gal g (y)$ for all
$y\in \ldom{B\gal  g}$. Since, $f\geq B\gal g$, 
then $f(y)=B\gal g(y)=+\infty$ for all $y\in Y\setminus \ldom{B\gal  g}$.
Therefore, $f=B\gal g$.
\end{proof}

We now prove Theorem~\ref{theo-uniquenew}.

\begin{proof}[Proof of Theorem~\ref{theo-uniquenew} in Cases~\eqref{lodiscrete}
and~\eqref{locompact}]
Assume that $Bf=g$ has a unique solution $f\in\sF$ and that
Condition~\eqref{lodiscrete} or Condition~\eqref{locompact} holds.
In particular, $B B\gal g=g$  and $X'=\dom{\partial\gal g}$. 
Set $Y'=\ldom {B\gal g}$.
By the equivalence \eqref{e-c0}$\Leftrightarrow$\eqref{e-c2}
in Theorem~\ref{cover}, 
$\{(\partial\gal g)^{-1}(y)\}_{y\in Y'}$ is a covering of $X'\cap \udom g$.
We shall show that this covering is topologically minimal.
Arguing by contradiction, and using the same arguments as in
the proof of the implication {\rm (\ref{fe-eq1},\ref{e-c2})}%
$\Rightarrow$\eqref{fe-c1} in Theorem~\ref{minicover}, we obtain that 
there exists an open set $U$ of $Y$ such that $U\cap Y'\neq\emptyset$ and
\eqref{eq-n} holds.
Since $Y$ is locally compact (this holds trivially in Case~\eqref{lodiscrete},
and by assumption in Case~\eqref{locompact}),
possibly after replacing $U$ by an open subset,
we can assume that $U$ is relatively compact, which means that its
closure is compact.
Taking $f=B\gal g$ on $Y\setminus U$ and $f=+\infty $ on $U$,
we deduce, as in the proof of {\rm (\ref{fe-eq1},\ref{e-c2})}%
$\Rightarrow$\eqref{fe-c1}, that $f\neq B\gal g$,
$Bf\leq g$ and $Bf=g$ on $X'$ (we obtain $Bf\leq g$ and
$Bf=g$ on $X'\cap \udom g$, and since
$Bf\leq g \Rightarrow Bf=g$ on $g^{-1}(-\infty )$,
we get $Bf=g$ on $X'$), hence
\begin{equation}
S\bydef\set{x\in X}{Bf(x)\neq g(x)}=\set{x\in X}{Bf(x)<g(x)}\subset
X\setminus X'\enspace.
\label{good}
\end{equation}

It remains to check that $S=\emptyset$, in order to 
contradict the uniqueness of the solution $f$  of $Bf=g$. 
Assume by contradiction that $S\neq \emptyset$.
If $x\in S$, it follows from~\eqref{keyprop}, that
$\sup_{y\in Y\setminus U} b(x,y, B\gal g(y))\leq Bf(x)<g(x)$.
Since $B B\gal g =g$, we get 
\[g(x)=\sup_{y\in U} b(x,y,B\gal g(y))\enspace.\]
Since $U$ has a compact closure
$\overline{U}$, and the function $y\mapsto 
b(x,y,B\gal g(y))$ is u.s.c.\ (this holds trivially in Case~\eqref{lodiscrete},
and this follows from Lemma~\ref{lem-compose} and
Assumption~\ref{a-bgalgisfinite}, in Case~\eqref{locompact}),
\begin{equation}
\exists y\in \overline{U} \text{ such that } g(x)=b(x,y,B\gal g(y))\enspace.
\label{eq0}\end{equation}
Since $g(x)>-\infty$, we get that  $(x,y)\in\sS$ and,
by~\eqref{e-inverse1} and Proposition~\ref{prop-partial},
 $x\in \partial B\gal g(y)=(\partial\gal g)^{-1}(y)$,
which implies
that $x\in\dom{\partial\gal g}=X'$. By \eqref{good}, we get 
a contradiction. 
\end{proof}

\begin{proof}[Proof of Theorem~\ref{theo-uniquenew}
in Case~\eqref{regular}]
Assume that $Bf=g$ has a unique solution $f\in\sF$ and that
Condition~\eqref{regular} holds.
In particular, $B B\gal g=g$.
Set $Y'=\ldom {B\gal g}=\dom {B\gal g}$ and $X'=\idom{g}$.
Theorem~\ref{cover} shows that 
$\{(\partial\gal g)^{-1}(y)\}_{y\in Y'}$ is a covering of 
$X'\cap \udom g=\idom{g}$. We shall show that this covering is minimal.

Arguing by contradiction, and using the same arguments as in
the proof of the implication {\rm (\ref{fe-eq1},\ref{e-c2})}%
$\Rightarrow$\eqref{fe-c1} in Theorem~\ref{minicover}, we obtain that there
exists an open set $U$ of $Y$ such that $U\cap Y'\neq\emptyset$ and
\eqref{eq-n} holds.
For any given basis of open neighbourhoods in $Y$, $\sB$,
possibly after replacing $U$ by an open subset,
we can assume that $U\in \sB$ and $U\cap Y'\neq \emptyset$.
We shall take either  $\sB$ as in condition \condc, or
$\sB$ as the basis of relatively compact open sets.

Fix $\varepsilon>0$ and consider the l.s.c.\ finite function
$w:Y\to [0,1]$ given by $w(y)=0$ for
$y\in Y\setminus U$ and $w(y)=\varepsilon$
for $y\in U$. Taking $f=B\gal g+ w$, 
we get that $f$ is l.s.c., $f$ satisfies~\eqref{keyprop},
and $f\neq B\gal g$ since
$f(y)=B\gal g(y)+\varepsilon>B\gal g(y)$ for $y\in
\dom{B\gal g}\cap U= Y'\cap U\neq\emptyset$.
As in the proof of the implication {\rm (\ref{fe-eq1},\ref{e-c2})}%
$\Rightarrow$\eqref{fe-c1}, we deduce that
$Bf\leq g$, $Bf=g$ on $X'$, hence \eqref{good} holds
and it remains to show that $S=\emptyset$.

Assume by contradiction that $S\neq \emptyset$.
We first prove that
\begin{subequations}\label{semi2}
\begin{align}
\label{semi2-1}\ldom g &\subset \ldom{Bf}\enspace, \\
\label{semi2-2}\udom g &\subset \udom{Bf}\enspace.\end{align}
\end{subequations}

Indeed,~\eqref{semi2-1} follows from $Bf \leq g$.
Let $x\in \udom g$, hence  $g(x)>-\infty$.
Since $g(x)=\sup_{y\in Y}
b(x,y,B\gal g(y))$, there exists $y\in Y$ such that 
$b(x,y,B\gal g(y))>-\infty$, hence $(x,y)\in \sS$ and $B\gal g(y)<+\infty$,
then $f(y)\leq B\gal g(y)+\varepsilon<+\infty$ and 
$Bf(x)\geq b(x,y,f(y))>-\infty$, which concludes the
proof of~\eqref{semi2-2}. 

Since $S\subset \ldom{Bf}\cap \udom g$, we deduce from~\eqref{semi2-2},
that $S\subset \dom{Bf}$.
Since $B$ is regular, $Bf$ and $g$ are l.s.c.\ on $X$ and quasi-continuous
on their domain. Hence, by Lemma~\ref{semiopen}, $S$ is semi-open relatively to
$\dom{Bf}$. Since $S\neq\emptyset$, $S$ has a nonempty interior relatively to
$\dom{Bf}$. This means that there exists an open set $V$ of $X$ such that
\begin{equation}
\label{semi1}
\emptyset\neq V\cap\dom{Bf}\subset S\enspace.\end{equation}

By \eqref{semi1} and~\eqref{semi2}, we get
\begin{equation}
\label{semi3}V\cap  \dom g\subset S\enspace.\end{equation}
If we know that $V\cap\dom g\neq\emptyset$, then since we assumed that
$\idom g$ is dense in $\dom g$, we get $V\cap\idom g\neq\emptyset$, so
by \eqref{semi3}, $S\cap \idom g\neq\emptyset$, i.e.\ $S\cap X'\neq \emptyset$,
which contradicts~\eqref{good}.
It remains to show that $V\cap \dom g\neq\emptyset$.

If $Y$ is locally compact, the arguments of the proof of
Theorem~\ref{theo-uniquenew} in Case~\eqref{locompact} show 
that \eqref{eq0} holds for all $x\in S$.
Since $B\gal g(y)>-\infty $ for all $y\in Y$, we deduce that
$g(x)<+\infty$, hence $S\subset \ldom g$.
Since we also have $S\subset \udom g$, we get $S\subset \dom g$
and by \eqref{semi1} and \eqref{semi2},
$V\cap \dom g=V\cap \dom{Bf}\neq\emptyset$.

Otherwise, $b$ satisfies Condition~\condc, and if $U\in\sB$ and $\varepsilon>0$
is chosen as in~\condc, we get that for all $x\in X$,
\begin{eqnarray*}
g(x)\geq Bf(x)&=&\sup_{y\in Y,\; B\gal g(y)<+\infty} b(x,y,B\gal g(y)+w(y))\\
&\geq& \sup_{y\in \sS_x,\; B\gal g(y)<+\infty} b(x,y,B\gal g(y))\\
&& \quad +\inf_{y\in  \sS_x,\; B\gal g(y)<+\infty}
 (b(x,y,B\gal g(y)+w(y))-b(x,y,B\gal g(y)))\\
&\geq& g(x)+\inf_{y\in  \sS_x\cap U,\; B\gal g(y) \in\R}
 (b(x,y,B\gal g(y)+\varepsilon)-b(x,y,B\gal g(y)))\\
&\geq& g(x)+\inf_{y\in  \sS_x\cap U,\; \alpha \in\R}
 (b(x,y,\alpha+\varepsilon)-b(x,y,\alpha))\enspace.
\end{eqnarray*}
By~\condc, we obtain $\dom g=\dom{Bf}$, which shows, by~\eqref{semi1},
that $V\cap \dom g\neq\emptyset$.
\end{proof}

\section{Algorithmic Issues}
\label{gal-ex-sec}
When $X$ and $Y$ are finite, and when the kernels
$b$ and $\bgal$ are given in an effective way,
Corollaries~\ref{cor-exis3} and~\ref{unifini} yield
an algorithm, which extends Zimmermann's algorithm,
to solve the equation $Bf=g$
and to decide the uniqueness of its solution.
Let us illustrate this algorithm by taking $X=\{x_1,x_2\}$
and $Y=\{y_1,y_2,y_3\}$.
Consider the kernel $b$ and the map $g$
given by the following table
\begin{align}
\label{e-image}
b:\;
\bordermatrix{   & y_1 & y_2 & y_3 \cr 
		 x_1 & -\lambda & 4-3\lambda & 2-\lambda\cr
		 x_2 & -\myb & 3-\lambda & -\lambda}\enspace,
\qquad 
g:\; \bordermatrix{ & \cr x_1 & 8\cr x_2 & 6}\enspace,
\end{align}
which means for instance that $b(x_1,y_2,\lambda)=4-3\lambda $
and $g(x_1)=8$. (We denote by $\sgn(\lambda)\in \{0,\pm 1\}$ the sign
of a scalar $\lambda$.)
Let $B:\RBY\to \RBX$ denote the (dual) functional Galois connection 
with kernel $b$.
Assumptions~{\rm (}\ref{a1}--\ref{assump4}{\rm )}
are clearly satisfied with $\sS=X\times Y$.
Then the kernel of $B\gal$ is
\[
\bgal:\;
\bordermatrix{   & x_1 & x_2 \cr 
		 y_1 & -\lambda & -\mybinv\cr
		 y_2 &  (4-\lambda)/3 & 3-\lambda \cr
		 y_3 & 2-\lambda & -\lambda }\enspace,
\]
and $B\gal g$ is given by:
\begin{align}
B\gal g: \;
\bordermatrix{   & x_1 & & x_2 \cr 
		 y_1 & -8 &\vee &  \underline{-\sqrt 6}\cr
		 y_2 &  \underline{(4-8)/3} & \vee &  3-6 \cr
		y_3 & \underline{2-8} & \vee &  \underline{-6} }
= \begin{pmatrix} -\sqrt 6 \\ -4/3 \\ -6 \end{pmatrix}\enspace,
\label{e-algo}
\end{align}
where we underlined the terms which determine the maximum
(recall that $\vee$ denotes the sup law).
By~\eqref{argmax2}, the sets 
$(\partial\gal g)^{-1}(y_j)$ can be read directly from~\eqref{e-algo}
by choosing, for each row $y_j$,
the $x_i$ variables corresponding to the underlined
terms:
\[
(\partial\gal g)^{-1}(y_1) =  \{x_2\}\enspace,\quad
(\partial\gal g)^{-1}(y_2) =  \{x_1\}\enspace,\quad
(\partial\gal g)^{-1}(y_3) =  \{x_1,x_2\}\enspace.
\]
Since the union of these subsets is equal to $X=\{x_1,x_2\}$,
it follows from Corollary~\ref{cor-exis3} that $f=B\gal g$ is
a solution of $Bf=g$. It follows from Corollary~\ref{unifini}
that this solution is not unique, because the covering
$\{ (\partial\gal g)^{-1}(y_j)\}_{1\leq j\leq 3}$ of $X$ 
is not minimal: for instance, 
$\{ (\partial\gal g)^{-1}(y_3)\}$
is a subcovering of $X$, which reflects the fact
that setting
$f(y_1)=f(y_2)=+\infty$ and $f(y_3)= -6$ yields
another solution of $Bf=g$.

More generally, a minimal covering
of a set of cardinality $n$ must consist of at most
$n$ sets, which implies that when $X$ and $Y$ are finite,
the number of elements of $Y$, i.e.\ the number of ``scalar unknowns'',
must not exceed the number
of elements of $X$, i.e.\ the number of ``scalar equations'',
for the solution of $Bf=g$ to be unique.

To show a uniqueness case,
consider the restriction $B_{1,2}:\RBYP\to \RBX$,
which is obtained by specialising
$B$ to those $f$ such that $f(y_3)=+\infty$.
Then the covering 
$\{(\partial\gal g)^{-1}(y_j)\}_{1\leq j \leq 2}$ of $X$ 
is minimal, which shows that setting $f(y_1)=-\sqrt 6$, $f(y_2)=-4/3$
yields the only solution of $B_{1,2}f=g$.

To illustrate the case where $Bf=g$ has no solution,
consider:
\[ 
g': \bordermatrix{ & \cr x_1 & 3\cr x_2 & -3} ,
\quad \mrm{ with } \quad 
B\gal g': \;
\bordermatrix{   & x_1 & & x_2 \cr 
		 y_1 & -3 &\vee &  \underline{\sqrt 3}\cr
		 y_2 &  (4-3)/3 & \vee &  \underline{3+3} \cr
		y_3 & 2-3 & \vee &  \underline{3} }
= \begin{pmatrix} \sqrt 3 \\ 6 \\ 3 \end{pmatrix}\enspace.
\] 
\sloppy
We see from Corollary~\ref{cor-exis3} that $Bf=g'$ has
no solution, because $\bigcup_{1\leq j\leq 3} 
(\partial\gal g')^{-1}(y_j)=\{x_2\}$
is not a covering of $X$. 

\begin{figure}[htb]
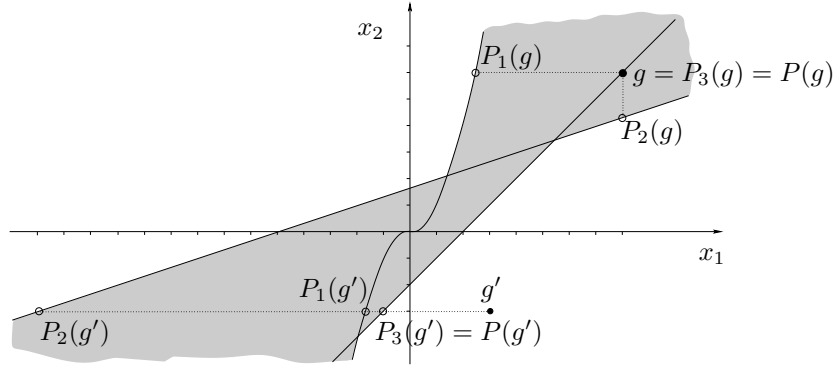

\input fig1
\caption{Image of the functional Galois connection~\eqref{e-image}.}
\label{fig1}
\end{figure}

Finally, let us interpret these computations
in geometric terms. For each $1\leq j\leq 3$, denote by $B_j$ the restriction
of $B$, $\RBYPj\to \RBX$, which is obtained by specialising
$B$ to those $f$ such that $f(y_k)=+\infty$
for $k\neq j$, so that the corresponding kernels
$b_j$ are given by:
\[
b_1:\; 
\bordermatrix{   & y_1 \cr 
		 x_1 & -\lambda \cr
		 x_2 & -\myb }\enspace,
\quad
b_2:\; 
\bordermatrix{   & y_2 \cr 
		 x_1 & 4-3\lambda\cr
		 x_2 & 3-\lambda}\enspace,
\quad
b_3:\; 
\bordermatrix{   & y_3 \cr 
		 x_1  & 2-\lambda\cr
		 x_2  & -\lambda}\enspace.
\]
The (set of finite points of the) image of the operator
$B_1$ is the curve $\lambda \to \sgn(\lambda)\lambda^2$
which is depicted on Figure~\ref{fig1}. The image
of $B_2$ (resp.\ $B_3$) is the line with slope $1/3$ (resp.\ $1$)
on the figure. The image of $B$ can be computed
readily from the images of $B_j$: since
$Bf=B_1f(y_1)\vee B_2f(y_2)\vee B_3f(y_3)$,
the image of $B$ is the sup-subsemilattice
of $\RBX$ generated by the images of $B_1,B_2,B_3$,
which corresponds to the gray region
on Figure~\ref{fig1}. 

Now, for each dual functional Galois connection $(B,B\gal)$, observe
that $BB\gal g$ is the maximum element of the image of $B$ which
is below $g$. Thus, $P=BB\gal$ is a nonlinear
projector on the image of $B$, and for $j=1,2, 3$,
consider the nonlinear projector $P_j =B_j B_j\gal$
on the image of $B_j$. By definition of Galois connections,
$P_j g(x)=b(x,y_j,B\gal g(y_j))$, thus
\(
P= \sup_{1\leq j\leq 3} P_j
\).
The element $g$, and its image by the projectors
$P_j$, are shown on Figure~\ref{fig1} (and can be computed directly
from the figure). For each $1\leq j\leq 3$,
the set $(\partial\gal g)^{-1}(y_j)$ represents
the subset of elements $x_i$ of $\{x_1,x_2\}$ 
such that $(P_j g)(x_i)=g(x_i)$. Thus,
the covering condition $X\subset \bigcup_{1\leq j \leq 3}
(\partial\gal g)^{-1}(y_j)$ is nothing but a combinatorial
rephrasing of $g=\sup_{1\leq j\leq 3} P_j g$.

\section{Some Examples of Moreau Conjugacies}\label{sec4}
We give now some applications of the results of Sections~\ref{exis-sec} 
and~\ref{uni-sec} to the case of the Moreau conjugacies $B$ and $B\gal$
given by~{\rm (\ref{e-repBB},\ref{defib})}, for
a kernel $\bar{b}$ taking only finite values. 
Since $\sS=X\times Y$, we have
 $B\gal g(y)>-\infty$ for all $y\in Y$
and $g\in\sG$ such that $g\not\equiv +\infty$.
\subsection{The Legendre-Fenchel Transform}\label{fenc-sec}
Let us consider the case where $X=Y=\R^n$ and $B=B\gal$ is the
Legendre-Fenchel transform, that is $B$ and $b$ are
given by~{\rm (\ref{e-repBB},\ref{defib})}, with
$\bar{b}(x,y)=\<x,y>$.
We have already shown in Section~\ref{sec31} that $b$ is continuous
in the second variable and coercive.
We also have:

\begin{lemma}\label{fenchel-reg}
The Legendre-Fenchel transform on $\R^n$ is regular.
\end{lemma}
\begin{proof} We need to show that for any function $f:\R^n\to\rbar$,
$g=f^\star$ is l.s.c.\ on $\R^n$ and is quasi-continuous on
its domain $\dom g$. We know that $g$ is either $\equiv +\infty$,
or $\equiv -\infty$ or a l.s.c.\ proper convex function.
Hence, it is l.s.c.\ and in the first two cases, the domain of $g$ 
is empty.
In the last case, since $g:\R^n\to \R\cup\{+\infty\}$ is l.s.c.\ 
on $\R^n$, and a fortiori on $\dom{g}$,
it is sufficient to prove that $g=\lsc (\usc (g))$
where $\lsc$ and $\usc$ envelopes are applied to the restrictions
to $\dom{g}$. Moreover, 
since $g$ is l.s.c., we get $g\leq \lsc (\usc (g))$, hence
it is sufficient to prove that $g\geq\lsc (\usc (g))$.
The following properties of l.s.c.\ proper convex functions can be
found in~\cite{ROCK}:
$g$ is continuous in the relative interior
of $\dom g$, that we denote by $\ridom{g}$ (recall that the relative interior
of a convex set is the interior of the
set for the topology of the affine hull of the set),
for any affine line $L$, the restriction of $g$ to $L$
is continuous on its domain $\dom{g}\cap L$, and
$\ridom{g}\cap L$ is dense in $\dom{g}\cap L$.
{From} this, we get that $\usc (g)=g$ on $\ridom{g}$. Let us
fix $x_0\in\ridom{g}$. For all $x\in\dom{g}$, take the affine line $L$
containing $x_0$ and $x$. Since there exists a sequence $(x_n)_{n\in\N}$
in $\ridom{g}\cap L$ converging to $x$ and since $g$ is continuous on
$\dom{g}\cap L$, we obtain that
$g(x)=\lim_{n\to \infty} g(x_n)= \lim_{n\to \infty} \usc(g)(x_n)
\geq \lsc (\usc (g))(x)$, which finishes
the proof.
\end{proof}

Using Theorems~\ref{cover} and~\ref{minicover1}, we get
\begin{proposition}
\sloppy
Let $g$ be a l.s.c.\ proper convex function on $\R^n$.
Then $\{ (\partial g)^{-1}(y)\}_{y\in \dom{g^\star}}$ is a covering
of  $\idom g$. 
Let $Z_a$ (resp.\ $Z_t$) be the set of algebraically (resp.\ topologically)
essential elements with respect to this covering,
and let $\intt{Z_t}$ denotes the interior of $Z_t$, relatively to
$\dom{g^\star}$.
Then \begin{align*}
& (f^{\star}\leq g\mrm{ and }f^{\star}(x)=g(x)\mrm{ for all }x\in\idom g)\\
& \qquad \implies (g^\star \leq f\mrm{ and } f(y)=g^\star(y)\mrm{ for all }
y\in Z_a\cup\intt{Z_t})\enspace.
\end{align*}
\end{proposition}

Since $Y$ is locally compact, and for any l.s.c.\ proper convex function
on $\R^n$ such that $\dom g$ has a nonempty interior, $\idom g$,
the set 
$\dom g$ is included in the closure of $\idom g$~\cite[Theorem 6.3]{ROCK},
we can also apply Theorem~\ref{theo-uniquenew}. We deduce:
\begin{proposition}\label{unic-fenchel1}
Let $g$ be a l.s.c.\ proper convex function on $\R^n$ such that
$\idom{g}\neq\emptyset$. The following statements are equivalent:
\begin{align}
& (f^{\star}\leq g\mrm{ and }f^{\star}(x)=g(x)\mrm{ for all }x\in\idom g)
\implies f=g^\star\enspace; \label{uf-1}\\
&\{ (\partial g)^{-1}(y)\}_{y\in \dom{g^\star} }
\mrm{ is a topologically minimal covering of } \idom g\enspace;
\label{uf-2}\\
&f^{\star}= g \implies f=g^\star\enspace.
\end{align}
\end{proposition}

The following classical notion is intermediate
between algebraic and topological minimality: 
a l.s.c.\ proper convex function $g$ on $\R^n$
is \NEW{essentially smooth}
if $\idom g\neq \emptyset$,
$g$ is differentiable in  $\idom g$, and 
the norm of the differential of $g$ at $x$ tends  to infinity, when
$x$ goes to the boundary of $\dom g$, see~\cite[Section~26]{ROCK}.
A l.s.c.\ proper  convex function $f$  on $\R^n$ is 
\NEW{essentially strictly convex}
if the restriction of  $f$ to any affine line (or segment) 
in $\dom f$ is strictly convex.
A l.s.c.\ proper convex function $g$ is essentially smooth if, and only if, 
its conjugate $g^{\star}$ is essentially strictly convex~\cite[Theorem 26.3]{ROCK}.
The following result, which can be compared 
with~\cite[Corollary~26.4.1]{ROCK},
is a corollary of Proposition~\ref{unic-fenchel1}.
It is the underlying argument of the G\"artner-Ellis theorem and
it is explicitly used in~\cite[Theorem 4.1 (c)]{VERV95}, 
in~\cite[Theorems 4.7 and 5.3]{gulinski}, and in~\cite[Lemmas
 3.2 and 3.5]{PUHAL94-1}.

\begin{corollary}
Let $g$ be an essentially smooth l.s.c.\ proper convex function on $\R^n$.
If $f$ is a l.s.c.\ function such that
$f^{\star}\leq g$ and $f^{\star}(x)=g(x)$ for all $x\in\idom g$, then $f=g^{\star}$.
In particular, $g$ has a unique preimage by the Legendre-Fenchel transform.
\end{corollary}
\begin{proof}
First, since $f^{\star}\leq g$ implies $f\geq g^{\star}$, so $f=g^{\star}$ outside
$\dom{g^{\star}}$, one can replace $Y$ by the affine hull
of $\dom{g^{\star}}$, so that $\idom{g^{\star}}$ is the relative interior of
$\dom{g^{\star}}$ (and is thus nonempty).
The conditions on the differentials of $g$ imply that $\partial g(x)$
is a singleton when $x\in \idom g$ and is empty elsewhere (see 
\cite[Theorem~26.1]{ROCK}).
Hence, applying Theorem~\ref{exis} to $g$,
we get that for all $y\in  \idom{g^{\star}}$ there exists $x\in\ldom g
=\dom g$
such that $y\in \partial g(x)$.
Since $\partial g(x)\neq\emptyset$, we get that
$x\in \idom g$, and $\{y\}= \partial g(x)$, hence any 
$y\in  \idom{g^{\star}}$ cannot be removed in the covering of $\idom g$ by 
$\{(\partial g)^{-1}(y)\}_{y\in \dom{g^{\star}}}$.
Moreover, since  $\dom{g^{\star}}$ is included in the closure
of $\idom{g^{\star}}$, any open set $U$ of  $\dom{g^{\star}}$ contains
a point $y\in \idom{g^{\star}}$, so $U$ cannot be removed
in the covering of $\idom g$ by 
$\{(\partial g)^{-1}(y)\}_{y\in \dom{g^{\star}}}$.
This shows that the covering of $\idom g$ by 
$\{(\partial g)^{-1}(y)\}_{y\in \dom{g^{\star}}}$
is topologically minimal.
The implication~\eqref{uf-2}$\Rightarrow$\eqref{uf-1} in 
Proposition~\ref{unic-fenchel1} yields the result of the corollary.
\end{proof}
\begin{example}
The following function $g$ satisfies~\eqref{uf-2} and thus~\eqref{uf-1},
but is not essentially smooth:
consider $X=Y=\R^2$,  $g=f^\star$ where $f:\R^2\to\R\cup\{+\infty\}$, 
with  $f(y)=y_1^2(y_2^2+3)$ if
$|y_2|\leq 1$ and  $f(y)=+\infty$ elsewhere.
Indeed, since $f$ is l.s.c.\ and convex, $f=g^{\star}$.
Since $f$ is not strictly convex on $y_1=0$, $g$ is not essentially
smooth. If $f$ is essentially strictly convex in a neighbourhood of $y\in \dom f$,
the point $y$ cannot be removed in the covering of $\dom{\partial g}$ by 
$\{(\partial g)^{-1}(y)\}_{y\in \dom{f}}$.
Then, since $\idom g=\dom{\partial g}=\dom g=\R^2$ and 
the loss of strict convexity of $f$ occurs only on a
line, any open set of $\dom f$ 
intersects the ``part'' of $\idom f$ where $f$ is essentially
strictly convex. This implies that the covering of $\idom g=\R^2$ by 
$\{(\partial g)^{-1}(y)\}_{y\in \dom{f}}$ is topologically minimal.
\end{example}

\subsection{Quadratic Kernels}\label{quadra-sect}
Let us consider the case where $Y=X=\R^n$ and
$b$ is given by~\eqref{defib} with $\bar{b}(x,y)=
b_a(x,y):= \<x,y>-\frac{a}{2} \|y\|^2$, where $\|\cdot\|$ is the Euclidean 
norm and $a\in \R$ is some constant.
Denoting $B_a$ and $B_a\gal$ the corresponding Moreau conjugacies 
given by~\eqref{e-repBB},
we get that $B_a f=(f+\frac{a}{2} \|\cdot\|^2)^\star$ and
 $B_a\gal g=-\frac{a}{2} \|\cdot\|^2+g^\star$, 
hence the properties of $B_a$ can be deduced
from those of the Legendre-Fenchel transform.
In particular:
\begin{corollary}\label{unic-fenchel2}
Let $g:\R^n\to\R\cup\{+\infty\}$ be 
an essentially smooth l.s.c.\ proper convex function.
If $f$ is a l.s.c.\ function such that
$B_a f\leq g$ and $B_a f(x)=g(x)$ for all $x\in\idom g$, then $f=B_a\gal g$.
\end{corollary}
Such kernels are useful, for instance, if we want to identify a function $f$
which is semiconvex but not convex. Indeed, if $f$ is semiconvex,
there exists $a\in \R$ such that $f+\frac{a}{2}\|\cdot\|^2$ is
strictly convex.
Hence, $g=B_a f$ satisfies the assumptions of Corollary~\ref{unic-fenchel2}.
Note that, by standard results of convex analysis, we know that when
$f$ is l.s.c.\ proper and convex, $B_a f$ is the inf-convolution
of $f^\star$ and $(\frac{a}{2} \|\cdot\|^2)^\star=\frac{1}{2a} \|\cdot\|^2$,
that is it is the Moreau-Yoshida regularisation of $f^\star$, which explain
why $B_a f$ is essentially smooth.

Similar results can be obtained when replacing the kernel $b_a$ by
$b'_a(x,y)= -\frac{a}{2} \|x-y\|^2$, with $a\neq 0$.

\subsection{$\omega$--Lipschitz Continuous Maps}\label{lip-sec}
Let $E$ be a Hausdorff topological vector space and
$\omega:E\to\R^+$ be a continuous subadditive map:
\[
\omega(x+y)\leq \omega(x)+\omega(y)\quad\mrm{for all } x,y\in E\enspace,
\]
such that $\omega(-x)=\omega(x)$ for
all $x\in E$, and $\omega(x)=0\Leftrightarrow x=0$. 
We say that a function $f:E\to \R$ is \NEW{$\omega$-Lipschitz 
continuous} if:
\[ |f(y)-f(x)|\leq \omega(y-x)\quad\mrm{ for all } x,y\in E\enspace,\]
and we denote by $\lipo{E}$ the set of $\omega$-Lipschitz 
continuous functions $f:E\to \R$.
If $E$ is a normed vector space with norm $\|\cdot\|$,
then $\omega(x)= a \|x\|^p$ satisfies the above properties for all
$a>0$ and $p\in (0,1]$, and in that case $\lipo{E}$ is the set of
H\"older continuous functions $f:E\to \R$ with exponent $p$ and
multiplicative constant $a$.

Take $Y=X=E$ and consider the kernel
$b$ given by~\eqref{defib} with $\bar{b}(x,y)=
b_\omega(x,y):= -\omega(y-x)$.
We denote by  $B_\omega $ and $B_\omega \gal$ the corresponding 
Moreau conjugacies given by~\eqref{e-repBB}.
We have $B_\omega =B_\omega \gal$.
When $\omega(x)= a \|x\|^p$, these Moreau conjugacies were studied by
Dolecki and Kurcyusz~\cite[Section 5]{dolecki}.
The kernel $b$ is continuous (in the second variable), but 
$b$ is not coercive, in general, since when $\omega =a \|\cdot\|$ and
$V$ is the ball of centre $x$ and radius $\varepsilon$,
$ b^\alpha_{x,V} (y)=\sup_{z\in V} a(\|y-x\|-\|y-z\|)+\alpha=a \varepsilon
+\alpha$.
We have:
\begin{proposition}\label{lip1}
Let $g\in\sG$. Then $g=B_\omega  B_\omega \gal g$ if, and only if,
either $g\equiv +\infty$, or $g\equiv -\infty$, or 
$g\in\lipo{E}$. In that case, we have $B_\omega \gal g=-g$.
\end{proposition}
\begin{proof}
Let $f\in\sF$. If $B_\omega f\not\equiv +\infty$
and $B_\omega f\not\equiv -\infty$, then $f(y)>-\infty$ for all
$y\in Y$ and there exists $y\in Y$ such that $f(y)<+\infty$.
Moreover, since $\omega$ is subadditive,
$b_\omega (\cdot,y)-\alpha $ is a $\omega$--Lipschitz continuous function
for all $y\in E$ and $\alpha\in \R$. Hence,
$B_\omega  f$ which is the supremum of a non-empty family
of $\omega$--Lipschitz continuous functions,
is also $\omega$-Lipschitz continuous,
which shows the ``only if'' part of the proposition.

Conversely, if $g\equiv +\infty$ or $g\equiv -\infty$,
then $B_\omega \gal g=-g$ and $g=B_\omega  B_\omega \gal g$.
Let $g\in\lipo{E}$.
Since $g(y)-g(x)\leq \omega(y-x)$ for all $x,y\in E$, we deduce:
\[ B_\omega \gal g(y)  =\sup_{x\in E} -\omega(y-x) -g(x)\leq -g(y)\enspace.\]
Moreover, taking $x=y$ in the supremum, we get that
$B_\omega \gal g(y)\geq -g(y)$, hence $B_\omega \gal g=-g$.
Since $B_\omega =B_\omega \gal$ and since $-g$ is also in $\lipo{E}$,
it follows, by application of the same argument,
that $B_\omega  (-g)= g$, hence
$B_\omega  B_\omega \gal g= B_\omega  (-g)=g$, which shows the ``if'' part 
of the proposition.
\end{proof}

Proposition~\ref{lip1} implies that $B_\omega =B_\omega \gal$ is regular, since any map
of the form $B_\omega  f$ is continuous. Also, 
$\dom{B_\omega  f}=\idom{B_\omega  f}$ is equal to $E$ or $\emptyset$
for all  $f\in\sF$. We also have:
\begin{lemma}\label{relcomp}
A map $f\in\sF$ is in $\sFC$ if, and only if,
$f+\omega$ has relatively compact finite sublevel sets,
which means that $\set{y\in E}{f(y)+\omega(y)\leq \beta}$ is
relatively compact, for all $\beta\in\R$.
\end{lemma}
\begin{proof}
By definition, $f\in\sFC$ if, and only if, 
$b_\omega (x,\cdot)-f$ has relatively compact finite superlevel sets.
Since $b_\omega (0,\cdot)=-\omega$, $f\in\sFC$ implies
that $f+\omega$ has relatively compact finite sublevel sets.
Conversely, if
$f+\omega$ has relatively compact finite sublevel sets,
then, by subadditivity of $\omega$,
$\set{y\in E}{b_\omega (x,y)-f(y)\geq \beta}\subset
\set{y\in E}{f(y)+w(y)\leq \omega(x)-\beta}$ is
relatively compact for all $\beta\in\R$ and $x\in E$.
\end{proof}

The condition of Lemma~\ref{relcomp} holds in particular when
$E=\R^n$, $\omega=a\|\cdot\|$ for some norm on $E$
and $f$ is Lipschitz continuous with a Lipschitz constant $b<a$.
We can apply Corollaries~\ref{cor-compact} and~\ref{uniccompact},
and Theorem~\ref{minicover1} 
to any  map $g$ such that $f=B_\omega \gal g$ satisfies the  
condition of Lemma~\ref{relcomp}.
Using Proposition~\ref{lip1}, we obtain:

\begin{corollary}
Let $g\in\lipo{E}$ be such that $g-\omega$
 has relatively compact finite superlevel sets,
then $\{ (\partial\gal g)^{-1}(y)\}_{y\in E}$ is a covering of  $E$. 
If 
\begin{equation}
\label{lipstrict}
|g(y)-g(x)|<\omega(y-x) \quad\mrm{for all } x,y\in E\mrm{ such that }
x\neq y\enspace,
\end{equation}
then $(f\in\sF\mrm{ and } B_\omega  f=g)\implies f=-g$.
\end{corollary}
\begin{proof}
The first assertion follows from Corollary~\ref{cor-compact}
and Lemma~\ref{relcomp}.
We shall prove that under~\eqref{lipstrict}, the covering is algebraically
(hence topologically) minimal, which will imply the last assertion of the
corollary, by using Corollary~\ref{uniccompact}.
Indeed, if $x,y\in E$, then 
\begin{equation}\label{lip2}
x\not\in \bigcup_{z\in E\setminus\{y\}} (\partial\gal g)^{-1}(z)
\iff  z\not \in \partial\gal g(x)\;\forall z\in E
\setminus\{y\} \iff \partial\gal g(x)\subset \{y\}\enspace.
\end{equation}
Since $g\in\lipo{E}$, we get $B_\omega \gal g=-g$, and, by~\eqref{e-inverse2},
\begin{align*}
\partial\gal g(x) &=\set{z\in E}{ B_\omega \gal g(z)= b_\omega (x,z)-g(x)}\\
&= \set{z\in E}{g(z)-g(x)=\omega(z-x)}\\
&= \set{z\in E}{g(z)-g(x)\geq \omega(z-x)}\enspace.
\end{align*}
Then $x\in \partial\gal g(x)$ for all $x\in E$, and
\begin{equation}\label{lip3}
\partial\gal g(x)\subset \{y\}\iff
x=y \mrm{ and }g(z)-g(y)<\omega(z-y)\;\forall z\in E\setminus\{y\}\enspace.
\end{equation}
Using~\eqref{lip2}, \eqref{lip3} and Definition~\ref{defi-cover2}, 
we get that the covering 
$\{ (\partial\gal g)^{-1}(y)\}_{y\in E}$ of $E$ is
algebraically  minimal if, and only if, for all $y\in E$, 
$g(z)-g(y)<\omega(z-y)$ for all $z\in E\setminus\{y\}$, which (by symmetry) is
equivalent to Condition~\eqref{lipstrict}.
\end{proof}

\subsection{L.s.c.\ Maps bounded from below}\label{lsc-sec}
Let $(E,\|\cdot\|)$ be a normed space, fix a constant $p>0$,
take $Y=E$, $X=E\times (0,+\infty)$ and consider the kernel
$b$ given by~\eqref{defib} with $\bar{b}(x,y)= 
 -x'' \|y-x'\|^p$, with $x=(x',x'')$,
$x'\in E$, and $x''\in (0,+\infty)$.
We denote by  $B$ and $B\gal$ the corresponding Moreau conjugacies
given by~\eqref{e-repBB}. 
These Moreau conjugacies were studied in~\cite[Section 4]{dolecki}.
When $p=1$, $B$ is used in~\cite{samb02} to define a (semi-)distance 
on the set of quasi-continuous functions from $E$ to $\R$.
When $p\leq 1$, $f\in\sF$, and $x=(x',x'')\in E\times (0,+\infty)$,
$Bf(x)=B_{\omega} f (x')$ 
where $B_{\omega}$ is given as in Section~\ref{lip-sec}
with $\omega= x''\|\cdot\|^p$, and
the results of this latter section show that $B$ is injective on
the set of H\"older continuous functions with exponent $p$.
When $E=\R^n$, $\|\cdot\|$ is the Euclidean norm, $p=2$, $f\in\sF$, and
$x=(x',x'')\in E\times (0,+\infty)$,
$Bf(x)=-x''\|x'\|^2+B_{2x''} f (2 x'' x')$ 
where $B_{2 x''}$ is given as in Section~\ref{quadra-sect}.
The results of this latter section show that $B$
is injective on the set of semiconvex maps.
Proposition~\ref{unicsci} below shows that indeed, for all $p>0$, 
$B$ is injective on a large set of l.s.c.\ functions. 
We first prove some preliminary results.

\begin{lemma}\label{lem-sci}
Let $f\in\sF$. Then either $Bf\equiv +\infty$, or $Bf\equiv -\infty$,
or there exists $a\geq 0$ such that 
\begin{equation}\label{dombf}
 E\times (a,+\infty)=\idom{Bf} \subset \dom{Bf}\subset
E\times [a,+\infty)\enspace.\end{equation}
More precisely, for all $x=(x',x''),\; z=(z',z'')\in X$ such that $x''>z''$,
we have
\begin{subequations}\label{lem-sci0}
\begin{align}
Bf(x)&\leq Bf(z)+K(x,z)\enspace, 
\\
\mrm{ with }  K(x,z) &:= \begin{cases}
 \left((z'')^{\frac{1}{1-p}}-(x'')^{\frac{1}{1-p}}
\right)^{1-p} \|x'-z'\|^p& \mrm{ when } p>1\enspace,\\
z''\|x'-z'\|^p& \mrm{ when } p\leq 1\enspace.
\end{cases}\label{lem-sci2}
\end{align}
\end{subequations}
\end{lemma}
\begin{proof}
Assume that $Bf\not\equiv +\infty$, and $Bf\not\equiv -\infty$.
Then there exists $y\in Y$ such that $f(y)<+\infty$,
which implies that $Bf(x)>-\infty$ for all $x\in X$,
hence $\dom{Bf}\neq \emptyset$.
Assume first that~\eqref{lem-sci0} is proved.
Let 
\[a=\inf\set{x''\in (0,+\infty) }{\exists x'\in E\mrm{ such that }
Bf(x',x'')<+\infty}\enspace.\]
Then $\dom{Bf}\subset E\times [a,+\infty)$, and since
$\idom{Bf}$ is open and included in $\dom{Bf}$,
we get that  $\idom{Bf}\subset E\times (a,+\infty)$.
Conversely, let $x_0=(x'_0,x''_0)\in E\times (a,+\infty)$. 
By definition of $a$, there exists
$z=(z',z'')\in  E\times (a,x''_0)$ such that $Bf(z)<+\infty$.
Let $\varepsilon=\frac{x''_0-z''}{2}>0$ and consider the
neighbourhood $V$ of $x_0$ given by
$V= \set{(x',x'')\in X}{\|x'-x'_0\|\leq \varepsilon\mrm{ and }
|x''-x''_0|\leq \varepsilon}$. Then $x''-z''\geq \varepsilon$ and
$\|x'-z'\|\leq \|x'_0-z'\|+\varepsilon$, for all $(x',x'')\in V$.
Using~\eqref{lem-sci0}, we obtain:
\begin{align*}
 \sup_{x\in V} Bf(x) &\leq Bf(z)+ \sup_{x\in V} K(x,z)\enspace,
\end{align*}
and 
\begin{align*}
\sup_{x\in V} K(x,z) &
 \leq \begin{cases}
 \left((z'')^{\frac{1}{1-p}}-(z''+\varepsilon)^{\frac{1}{1-p}}
\right)^{1-p} (\|x'_0-z'\|+\varepsilon)^p& \mrm{when } p>1\enspace,\\
z''(\|x'_0-z'\|+\varepsilon)^p& \mrm{when } p\leq 1\enspace,
\end{cases}
\end{align*}
hence $x_0\in \idom{Bf}$, which finishes the proof of~\eqref{dombf}.

Let us now prove~\eqref{lem-sci0}.
Let $x=(x',x''),\; z=(z',z'')\in X$ be such that $x''>z''$.
Using the definition of $Bf$, we deduce that
$Bf(x)\leq Bf(z)+K_0(x,z)$, where
\begin{eqnarray}
K_0(x,z) & =  &
\sup_{y\in E} \left( -x''\|y-x'\|^p+z''\|y-z'\|^p\right) \nonumber\\
&\leq &  \sup_{y\in E} \left( -x''\|y-x'\|^p+z''
(\|y-x'\|+\|x'-z'\|)^p\right) \nonumber\\
&\leq & \sup_{\rho\geq 0} \left( -x''\rho^p+z''
(\rho+\|x'-z'\|)^p\right)\enspace.\label{majp}
\end{eqnarray}
Computing the supremum in~\eqref{majp}, we obtain that 
$K_0(x,z)\leq K(x,z)$ with $K$ given by~\eqref{lem-sci2}, 
which shows~\eqref{lem-sci0}.
\end{proof}

\begin{proposition}\label{unicsci}
Let $f\in\sF$. Then $f=B\gal B f$ if, and only if,
$f\equiv -\infty$, or there exists $a>0$ such that 
$f+a \|\cdot\|^p$ is bounded from below.
Moreover, $\dom{Bf}\neq\emptyset\implies f=B\gal B f$.
\end{proposition}
\begin{proof}
If $f=B\gal B f$ and $f\not\equiv -\infty$, then there exists
$x\in X$ such that $Bf(x) <+\infty$. By Lemma~\ref{lem-sci},
either $Bf\equiv -\infty$, or~\eqref{dombf} holds.
In the first case, $f\equiv +\infty$, thus $f+a \|\cdot\|^p$ is bounded from below, for all $a>0$.
In the second case,  
there exists $a\geq 0$ such that $\idom{Bf}=E\times (a,+\infty)$.
Taking $a'>a$, we get that $Bf(0,a')<+\infty$, which means
that $f+a'\|\cdot\|^p$ is bounded from below.

Conversely, if $f\equiv -\infty$ or $f\equiv +\infty$, then
$f=B\gal B f$. Assume that $f\in\sF$ is such that
$f+a\|\cdot\|^p$ is bounded from below,
for some $a>0$, and that $f\not\equiv +\infty$.
Then, $Bf(0,a)<+\infty$ and $Bf(x)>-\infty$ for all $x\in X$,
which shows that $\dom{Bf}\neq\emptyset$.
It remains to prove the last assertion of the proposition,
that we shall derive from a result of~\cite{dolecki}.
Recall that if $\Phi$ is a set of real valued functions
on $E$, a function $f$ is $\Phi$-convex if it can
be written as a pointwise supremum of a possibly
infinite family of elements of $\Phi$.
Let us take for $\Phi$ the set
of functions from $E$ to $\R$ 
of the form $y\mapsto -x''\|y-x'\|^p+r$ with $x''>0$,
$r\in \R$ and $x'\in E$. When $\dom{Bf}\neq \emptyset$,
$f$ is bounded from below by an element of $\Phi$.
Then, by~\cite[Theorem 4.2]{dolecki}, $f$ is $\Phi$-convex,
which implies that $f=B\gal B f$.
\end{proof}

In order to deduce covering properties as in Section~\ref{uni-sec},
we need to show some properties of $b$ and $B$.
First, it is clear that $b$ is continuous (in the second variable).
The following result may be compared with~\cite[Lemma 4.4]{dolecki}
\begin{lemma}\label{lem-coerc}
When $E=\R^n$, the kernel $b$ is coercive.
\end{lemma}
\begin{proof}
Let $x\in E$, $\alpha\in\R$ and $V$ be a neighbourhood of $x$.
The map $b_{x,V}^\alpha$ defined in~\eqref{e-def-bv} satisfies:
\begin{align*}
 b_{x,V}^{\alpha}(y)& = 
\sup_{z\in V} \bar{b}(z,y)-\bar{b}(x,y)+\alpha 
=\sup_{(z',z'')\in V} -z'' \|y-z'\|^p+x'' \|y-x'\|^p+
\alpha\enspace.
\end{align*}
Let $\varepsilon >0$ be such that $V\supset  \set{(z',z'')\in X}{
\|z'-x'\|\leq \varepsilon,\; |z''-x''|\leq \varepsilon}$. We get that
$ b_{x,V}^{\alpha}(y)\geq \varepsilon \|y-x'\|^p+\alpha$,
hence $b_{x,V}^\alpha$ has bounded sublevel sets.
When $E=\R^n$, this implies that 
$b_{x,V}^\alpha$ has relatively compact  sublevel sets, hence
$b$ is coercive.
\end{proof}

\begin{lemma}\label{lem-regul}
$B$ is regular.
\end{lemma}
\begin{proof}
Let $f\in \sF$ and $g=Bf$. Then $g$ is l.s.c.\ as the supremum of 
continuous maps. To show that $g$ is quasi-continuous on its domain
it suffices to prove that $g=\lsc (\usc (g))$
where $\lsc$ and $\usc$ envelopes are applied to 
the restrictions to $\dom{g}$. Moreover, 
since $g$ is l.s.c., we get $g\leq \lsc (\usc (g))$, hence
it is sufficient to prove that $g\geq\lsc (\usc (g))$.
This is true if $g\equiv +\infty$ or $g\equiv -\infty$.
Otherwise,~\eqref{dombf} and~\eqref{lem-sci0} hold and
$\udom{g}=X$.
In particular, for all fixed $z=(z',z'')\in \dom g$,
$E\times (z'',+\infty)$ is open and included in  $\dom g$,
$g(x)\leq g(z)+K(x,z)$ for all $x\in E\times (z'',+\infty)$,
and $x\mapsto K(x,z)$ is continuous on $E\times (z'',+\infty)$.
Hence
\begin{equation}\label{reg-1}
\usc (g)(x)\leq g(z)+ K(x,z)\mrm{ for all }z=(z',z'')\in \dom g,\;
x\in E\times (z'',+\infty)\enspace.
\end{equation}
Since $\idom{g}$ is the interior of $\dom{g}$,
 for all $x=(x',x'') \in \idom{g}$, 
there exists $\varepsilon>0$ such that $z=(x',x''-\varepsilon)\in \dom{g}$.
Hence, using~\eqref{reg-1}, we get that 
$\usc(g)(x)\leq g(z)+K(x,z)=g(z)$.
It follows that 
\begin{equation}\label{reg-2}
\usc (g)(x)\leq \lim_{\varepsilon\to 0^+}  g(x',x''-\varepsilon)\enspace.
\end{equation} 
Moreover, it is clear that for all fixed $x'\in E$,
$x''\in (0,+\infty) \mapsto g(x',x'')$ is a nonincreasing 
l.s.c.\ proper convex map. Since $(0,+\infty)$ is one dimensional,
this implies in particular that this map is continuous on its 
domain (see~\cite{ROCK}).
Therefore, it follows from~\eqref{reg-2} that 
$\usc(g)(x)\leq g(x)$ for all $x\in\idom{g}$.
This shows that $g$ is continuous in the interior of its domain.
Now, by~\eqref{dombf}, if $x=(x',x'')\in\dom{g}$, 
$(x',x''+\varepsilon)\in\idom{g}$ 
for all $\varepsilon>0$, 
and since $x''\in (0,+\infty) 
\mapsto g(x',x'')$ is continuous on its  domain, we get 
\[ \lsc(\usc (g))(x)\leq \liminf_{\varepsilon\to 0^+} 
\usc(g)(x',x''+\varepsilon) =\liminf_{\varepsilon\to 0^+} 
g(x',x''+\varepsilon) = g(x)\enspace,\]
which finishes the proof of $\lsc(\usc (g))\leq g$.
\end{proof}
\begin{corollary}
Assume that $E=\R^n$. Let $g\in\sG$ such that $\dom{g}\neq\emptyset$
and $g=BB\gal g$.
Then $(f\in\sF\mrm{ and } Bf=g)\;\implies \; f=B\gal g$.
Moreover 
$\{ (\partial\gal g)^{-1}(y)\}_{y\in \dom{B\gal g}}$ is a topological minimal
covering of  $\idom g$. 
Finally, if $B\gal g$ is quasi-continuous on its domain, then 
$(Bf\leq g\mrm{ and } Bf=g\mrm{ on } \idom{g})\implies f=B\gal g$.
\end{corollary}
\begin{proof}
The first assertion of the corollary follows from the last assertion of 
Proposition~\ref{unicsci}, since $\dom{g}\neq\emptyset$.
This shows that Problem $(\sP)$ has a unique solution.
Since $b$ is continuous in the second variable and $\dom g\neq \emptyset$,
Assumption~\ref{a-bgalgisfinite} holds.
Moreover, $b$ is coercive (by Lemma~\ref{lem-coerc}),
$B$ is regular (by Lemma~\ref{lem-regul}), Condition~\condc\ holds,   
and $\dom{g}$ is included in the closure of $\idom{g}$
(by~\eqref{dombf}). Hence,
applying Theorem~\ref{theo-uniquenew} in Case~\eqref{regular},
we obtain the second assertion of the corollary.
Finally, applying the implication~\eqref{fe-c1}$\Rightarrow$\eqref{fe-eq1}
in Theorem~\ref{minicover} to $X'=\idom g$, we obtain
the last assertion of the corollary.
\end{proof}

\newcommand{\etalchar}[1]{$^{#1}$}
\providecommand{\bysame}{\leavevmode\hbox to3em{\hrulefill}\thinspace}


\begin{thebibliography}{GHK{\etalchar{+}}80}

\bibitem[AGK02]{AGK1}
M.~Akian, S.~Gaubert, and V.~Kolokoltsov, \emph{Invertibility of functional
  {G}alois connections}, C. R. Acad. Sci. Paris \textbf{Ser. {I} 335} (2002),
  1--6.

\bibitem[AGK04]{AGK2}
\bysame, \emph{Invertibility of {M}oreau conjugacies and large deviations}, In
  preparation, 2004.

\bibitem[Aki99]{DENSITE}
M.~Akian, \emph{Densities of idempotent measures and large deviations},
  Transactions of the American Mathematical Society \textbf{351} (1999),
  no.~11, 4515--4543.

\bibitem[Bal77]{balder}
E.~J. Balder, \emph{An extension of duality-stability relations to nonconvex
  optimization problems}, SIAM J. Control Optimization \textbf{15} (1977),
  no.~2, 329--343.

\bibitem[BCOQ92]{BCOQ}
F.~Baccelli, G.~Cohen, G.~J. Olsder, and J.-P. Quadrat, \emph{Synchronization
  and linearity~: an algebra for discrete events systems}, John Wiley \& Sons,
  New-York, 1992.

\bibitem[Bir95]{BIRK}
G.~Birkhoff, \emph{Lattice theory}, Colloquium publications, vol.~25, American
  Mathematical Society, Providence, 1995, (first edition, 1940).

\bibitem[BJ72]{BLY}
T.~S. Blyth and M.~F. Janowitz, \emph{Residuation theory}, Pergamon Press,
  1972.

\bibitem[But94]{butkovip94}
P.~Butkovi\v{c}, \emph{Strong regularity of matrices --- a survey of results},
  Discrete Applied Mathematics \textbf{48} (1994), 45--68.

\bibitem[But00]{butkovic}
P.~Butkovi\v{c}, \emph{Simple image set of (max, +) linear mappings}, Discrete
  Appl. Math. \textbf{105} (2000), no.~1-3, 73--86.

\bibitem[CG79]{cuning}
R.~A. Cuninghame-Green, \emph{Minimax algebra}, Lecture notes in Economics and
  Mathematical Systems, no. 166, Springer, 1979.

\bibitem[CGQ04]{cgq02}
G.~Cohen, S.~Gaubert, and J.P. Quadrat, \emph{Duality and separation theorem in
  idempotent semimodules}, Linear Algebra and Appl. \textbf{379} (2004),
  395--422.

\bibitem[DJLC53]{DUB}
M.-L. Dubreuil-Jacotin, L.~Lesieur, and R.~Croisot, \emph{Le\c{c}ons sur la
  th\'eorie des treillis des structures alg\'ebriques ordonn\'ees et des
  treillis g\'eom\'etriques}, Gauthier-Villars, Paris, 1953.

\bibitem[DK78]{dolecki}
S.~Dolecki and S.~Kurcyusz, \emph{On {$\Phi $}-convexity in extremal problems},
  SIAM J. Control Optimization \textbf{16} (1978), no.~2, 277--300.

\bibitem[DZ93]{DEMBO}
A.~Dembo and O.~Zeitouni, \emph{Large deviations techniques and applications},
  Jones and Barlett, Boston, MA, 1993.

\bibitem[Eve44]{everett}
C.~J. Everett, \emph{Closure operators and galois theory in lattices},
  Transactions of the American Mathematical Society \textbf{55} (1944),
  514--525.

\bibitem[GHK{\etalchar{+}}80]{LATT}
G.~Gierz, K.~H. Hoffman, K.~Keimel, J.~D. Lawson, M.~Mislove, and D.~S. Scott,
  \emph{A compendium of continuous lattices}, Springer Verlag, Berlin, 1980.

\bibitem[GM01]{gondran}
M.~Gondran and M.~Minoux, \emph{Graphes, dio\"\i des et semi-anneaux}, TEC \&
  DOC, Paris, 2001.

\bibitem[GR02]{goebel}
R.~Goebel and R.~T. Rockafellar, \emph{Generalized conjugacy in
  {H}amilton-{J}acobi theory for fully convex {L}agrangians}, Journal of Convex
  Analysis \textbf{9} (2002), no.~2, 463--473.

\bibitem[Gul03]{gulinski}
O.~V. Gulinsky, \emph{The principle of the largest terms and quantum large
  deviations}, Kybernetika \textbf{39} (2003), no.~2, 229--247.

\bibitem[KM87]{KOLO88}
V.~N. Kolokoltsov and V.~P. Maslov, \emph{The general form of the endomorphisms
  in the space of continuous functions with values in a numerical commutative
  semiring (with the operation $\oplus=\max$)}, Dokl. Akad. Nauk SSSR
  \textbf{295} (1987), no.~2, 283--287, Engl. transl. in Sov. Math. Dokl., 36
  (1), 55-59 (1988).

\bibitem[KM97]{maslovkolokoltsov95}
V.~Kolokoltsov and V.~Maslov, \emph{Idempotent analysis and applications},
  Kluwer Acad. Publisher, 1997.

\bibitem[Kol90]{KOLO90}
V.~N. Kolokoltsov, \emph{Maslov's arithmetic in general topology},
  In~V.~Fedorchuk et al (Eds.) Geometry, Topology and Applications, Moskov.
  Instrument. Inst., 1990, (in Russian), pp.~64--68.

\bibitem[Kol92]{kolokoltsov92}
V.~Kolokoltsov, \emph{On linear, additive, and homogeneous operators in
  idempotent analysis}, {in~\cite{maslov92}}, 1992, pp.~87--101.

\bibitem[Lin79]{lindberg}
P.~O. Lindberg, \emph{A generalization of {F}enchel conjugation giving
  generalized {L}agrangians and symmetric nonconvex duality}, Survey of
  mathematical programming (Proc. Ninth Internat. Math. Programming Sympos.,
  Budapest, 1976), Vol. 1, North-Holland, Amsterdam, 1979, pp.~249--267.

\bibitem[ML88]{martinez88}
J.-E. Mart{\'{\i}}nez-Legaz, \emph{Quasiconvex duality theory by generalized
  conjugation methods}, Optimization \textbf{19} (1988), no.~5, 603--652.

\bibitem[MLS90]{martinez90}
J.-E. Mart{\'{\i}}nez-Legaz and I.~Singer, \emph{Dualities between complete
  lattices}, Optimization \textbf{21} (1990), no.~4, 481--508.

\bibitem[MLS95]{martinez95}
\bysame, \emph{Subdifferentials with respect to dualities}, ZOR---Math. Methods
  Oper. Res. \textbf{42} (1995), no.~1, 109--125.

\bibitem[Mor70]{moreau70}
J.-J. Moreau, \emph{Inf-convolution, sous-additivit\'e, convexit\'e des
  fonctions num\'eriques}, J. Math. Pures Appl. (9) \textbf{49} (1970),
  109--154.

\bibitem[MS92]{maslov92}
V.~Maslov and S.~Samborski\u\i (eds.), \emph{Idempotent analysis}, Adv. in Sov.
  Math., vol.~13, AMS, RI, 1992.

\bibitem[Neu89]{NEUBRUNN}
T.~Neubrunn, \emph{Quasi-continuity}, Real Anal. Exchange \textbf{14}
  (1988/89), no.~2, 259--306.

\bibitem[Ore44]{ore}
O.~Ore, \emph{Galois connexions}, Transactions of the American Mathematical
  Society \textbf{55} (1944), 493--513.

\bibitem[OV95]{VERV95}
G.~L. O'Brien and W.~Vervaat, \emph{Compactness in the theory of large
  deviations}, Stochastic processes and their applications \textbf{57} (1995),
  1--10.

\bibitem[Puh94]{PUHAL94-1}
A.~Puhalskii, \emph{The method of stochastic exponentials for large
  deviations}, Stochastic processes and their applications \textbf{54} (1994),
  45--70.

\bibitem[Roc70]{ROCK}
R.~T. Rockafellar, \emph{Convex analysis}, Princeton University Press
  Princeton, N.J., 1970.

\bibitem[RR98]{rachev}
S.~T. Rachev and L.~R\"uschendorf, \emph{Mass transportation problems, volume
  {I}: theory}, Springer, 1998.

\bibitem[R{\"u}s91]{frechet}
L.~R{\"u}schendorf, \emph{Fr\'echet-bounds and their applications}, Advances in
  probability distributions with given marginals (Rome, 1990), Math. Appl.,
  vol.~67, Kluwer Acad. Publ., Dordrecht, 1991, pp.~151--187.

\bibitem[R{\"u}s95]{ruschendorf95}
\bysame, \emph{Optimal solutions of multivariate coupling problems}, Appl.
  Math. (Warsaw) \textbf{23} (1995), no.~3, 325--338.

\bibitem[RW98]{rockwets}
R.~T. Rockafellar and R.~J.-B. Wets, \emph{Variational analysis},
  Springer-Verlag, Berlin, 1998.

\bibitem[Sam02]{samb02}
S.~Samborski, \emph{A new function space and extension of partial differential
  operators in it}, Research Report 2002-13, Universit\'e de Caen, Laboratoire
  SDAD, 2002.

\bibitem[Sin86]{singer86}
I.~Singer, \emph{Some relations between dualities, polarities, coupling
  functionals, and conjugations}, J. Math. Anal. Appl. \textbf{115} (1986),
  no.~1, 1--22.

\bibitem[Sin97]{singer}
\bysame, \emph{Abstract convex analysis}, John Wiley \& Sons Inc., New York,
  1997.

\bibitem[Sin02]{singer2}
\bysame, \emph{Further applications of the additive min-type coupling
  function}, Optimization \textbf{51} (2002), no.~3, 471--485.

\bibitem[Vil03]{villani}
C.~Villani, \emph{Topics in optimal transportation}, Graduate Studies in
  Mathematics, vol.~58, American Mathematical Society, Providence, 2003.

\bibitem[Vol98]{volle}
M.~Volle, \emph{Duality for the level sum of quasiconvex functions and
  applications}, ESAIM: Control, Optimisation and calculus of variations
  \textbf{3} (1998), 329--343.

\bibitem[Vor67]{vorobyev67}
N.~N. Vorobyev, \emph{Extremal algebra of positive matrices}, Elektron.
  Informationsverarbeitung und Kybernetik \textbf{3} (1967), 39--71, (in
  Russian).

\bibitem[Zim76]{Zimmermann.K}
K.~Zimmermann, \emph{Extrem\'aln\'\i\ algebra}, Ekonomick\'y \`ustav \u CSAV,
  Praha, 1976, (in Czech).

\end{thebibliography}
\end{document}